\title{Spectral Measures for $G_2$ {II}: finite subgroups}
\author{
{\sc David E.\ Evans and Mathew Pugh}\\
 {\footnotesize School of Mathematics, Cardiff University,}\\  {\footnotesize Senghennydd Road, Cardiff CF24 4AG, Wales, U.K.}
}
\date{\today}

\documentclass[12pt]{article}

\usepackage{amsmath,amssymb,amsthm}
\usepackage{graphicx}
\usepackage{framed}
\usepackage{color}
\usepackage[usenames,dvipsnames]{xcolor}
\usepackage{multirow}
\usepackage{cancel}
\usepackage{array}

\theoremstyle{definition}
\newtheorem{Def}{Definition}[section]

\newtheorem{Thm}[Def]{Theorem}
\newtheorem{Rem}[Def]{Remark}

\newcommand\bbR{\mathbb{R}}

\newcommand\bbZ{\mathbb{Z}}

\newcommand\bbT{\mathbb{T}}
\newcommand\ben{\begin{enumerate}}
\newcommand\een{\end{enumerate}}

\begin{document}
\maketitle

\begin{abstract}
Joint spectral measures associated to the rank two Lie group $G_2$, including the representation graphs for the irreducible representations of $G_2$ and its maximal torus, nimrep graphs associated to the $G_2$ modular invariants have been studied. In this paper we study the joint spectral measures for the McKay graphs (or representation graphs) of finite subgroups of $G_2$. Using character theoretic methods we classify all non-conjugate embeddings of each subgroup into the fundamental representation of $G_2$ and present their McKay graphs, some of which are new.
\end{abstract}

{\footnotesize
\tableofcontents
}

\section{Introduction} \label{sect:intro}

Two fundamental algebraic structures in conformal quantum field theories are modular tensor categories and conformal nets of factors. The representation theory of the conformal nets of factors  in the algebraic formulation of a conformal field theory, is a modular tensor category describing the Verlinde ring of primary fields. The full conformal field theory is described by extensions of the conformal net, studied through subfactors. The study of modular tensor categories and subfactors and their classification leads to the quest to develop new invariants associated to these categories. Among the most basic invariant for a semisimple tensor category is the underlying fusion or Grothendieck ring, whilst for a subfactor basic invariants which are not complete are the Jones index, the principal graph and its dual \cite{goodman/de_la_harpe/jones:1989, evans/kawahigashi:1998}. At a more sophisticated level for a subfactor is the standard invariant (also known as the $\lambda$-lattice \cite{popa:1995}, or planar algebra \cite{jones:planar}). A modular tensor category is described by its braided structure from which one can derive  modular data (in particular its $S$ and $T$ matrices). In the cases of module categories for modular tensor categories (or equivalently braided subfactors in the conformal net picture), invariants include its nimrep (Non-negative Integer Matrix REPresentation of the fusion rules of the underlying modular tensor category) and the cell system built on this nimrep. Another invariant, considered in the work presented in this paper, is the spectral measure associated to the nimrep \cite{banica/bisch:2007}.

The representation ring of a finite group has a basis indexed by simple representations with multiplication given by tensor product.  Each representation acts by left multiplication and it is natural to ask what the eigenvalues (with multiplicity) and eigenvectors are for this linear operator. For compact groups a Hilbert space structure is needed on the representation ring and instead of eigenvalues with multiplicity one needs to consider the spectral measure.

Spectral measures associated to the compact Lie groups $SU(2)$ and $SU(3)$ and their maximal tori, nimrep graphs associated to the $SU(2)$ and $SU(3)$ modular invariants, and the McKay graphs (or representation graphs) for finite subgroups of $SU(2)$ and $SU(3)$ were studied in \cite{banica/bisch:2007, evans/pugh:2009v, evans/pugh:2010i}.
Spectral measures associated to the compact Lie group $Sp(2)$ are studied in \cite{evans/pugh:2012iii}, whilst spectral measures associated to other compact rank two Lie groups and their maximal tori are studied in \cite{evans/pugh:2012iv}.

The real rank two Lie group $G_2 = \mathrm{Aut}(\mathbb{O})$ is the automorphism group of the octonions $\mathbb{O}$. It is a simply connected and compact Lie group of dimension 14 and is the smallest of the exceptional Lie groups. It is isomorphic to the subgroup of $SO(7)$ that fixes any particular vector in its 8-dimensional real spinor representation.
Spectral measures for the representation graphs for the irreducible representations of $G_2$ and its maximal torus $\bbT^2$, and for nimrep graphs associated to the $G_2$ modular invariants, were studied in \cite{evans/pugh:2012i}.
The spectral measure for $G_2$ itself appears in \cite{uhlmann/meinel/wipf:2007} as the reduced Haar measure, and is useful for the study of how gluons transform under $G_2$ at finite temperature \cite{wellegehausen/wipf/wozar:2009}. These measures are also useful in other applications, e.g. the computation of glueball masses in Hamiltonian lattice gauge theories or in studies of the strong coupling limit of these theories \cite{kluberg-stern/morel/petersson:1983}.

In this paper we consider the McKay graphs for finite subgroups of $G_2$ and study their spectral measures.
In the context of finite subgroups of a compact group, the spectral measure simply encodes the character table of the group. However in the context of quantum subgroups discussed above, the spectral measure provides the appropriate generalisation of the character table. We trust that framing this work in terms of spectral measures will make the case of quantum groups more accessible for those who are more familiar with the classical group case.

The spectral measure of a self-adjoint operator $a$ in a unital $C^{\ast}$-algebra $A$ with state $\varphi$ is a compactly supported probability measure $\nu_a$ on the spectrum $\sigma(a) \subset \bbR$ of $a$, uniquely determined by its moments
$$\varphi(a^m) = \int_{\sigma(a)} x^m \mathrm{d}\nu_a (x),$$
for all non-negative integers $m$. We need to have a concept to encode the corresponding notions for commuting operators. This is familiar in classical and even free probability \cite{guionnet:2010} as joint laws.

The self-adjoint operators we consider here are the adjacency matrices of the McKay graphs for finite subgroups of $G_2$. The characters of the rank two Lie group $G_2$ are functions on $\bbT^2$, and it is convenient to write the spectral measures for these operators as measures $\varepsilon$ over the torus $\bbT^2$. However, $\bbT^2$ has dimension one greater than $\sigma(a) \subset \bbR$, so that there is an infinite family of pullback measures $\varepsilon$ over $\bbT^2$ for any spectral measure $\nu_a$. The details of the relation between the measures $\varepsilon$ and $\nu_a$ are given in Section \ref{sect:measures-different_domains}.

In order to remove this ambiguity, we also consider joint spectral measures, that is, measures over the joint spectrum $\sigma(a,b) \subset \sigma(a) \times \sigma(b) \subset \bbR^2$ of commuting self-adjoint operators $a$ and $b$. The abelian $C^{\ast}$-algebra $B$ generated by $a$, $b$ and the identity 1 is isomorphic to $C(X)$, where $X$ is the spectrum of $B$. Then the joint spectrum is defined as $\sigma(a,b) = \{ (a(x), b(x)) | \, x \in X \}$. In fact, one can identify the spectrum $X$ with its image $\sigma(a,b)$ in $\bbR^2$, since the map $x \mapsto (a(x), b(x))$ is continuous and injective, and hence a homomorphism since $X$ is compact \cite{takesaki:2002}.
In the case where the operators $a$, $b$ act on a finite-dimensional Hilbert space, this is the set of all pairs of real numbers $(\lambda_a,\lambda_b)$ for which there exists a vector $\phi$, $||\phi||=1$, such that $a\phi = \lambda_a \phi$, $b\phi = \lambda_b \phi$.
Then there exists a compactly supported probability measure $\widetilde{\nu}_{a,b}$ on $\sigma(a,b)$, which is uniquely determined by its cross moments
\begin{equation} \label{eqn:cross_moments_sa_operators}
\varphi(a^m b^n) = \int_{\sigma(a,b)} x^m y^n \mathrm{d}\widetilde{\nu}_{a,b} (x,y),
\end{equation}
for all non-negative integers $m$, $n$.
Such joint spectral measures specialize to the spectral measures $\nu_a$ (respectively $\nu_b$) by integrating over all $y$ for which $(\lambda_a,y) \in \sigma(a,b)$ (respectively all $x$ for which $(x,\lambda_b) \in \sigma(a,b)$).
As discussed in Section \ref{sect:measures-different_domains}, such a measure uniquely defines a measure over $\bbT^2$ invariant under an action of the Weyl group of $G_2$. In this paper we determine the joint spectral measure for each finite subgroup $\Gamma$ of $G_2$ for all non-conjugate embeddings of $\Gamma$ into the fundamental representations of $G_2$.
In the process, we correct a number of errors that have previously appeared in the literature.

A number of recent papers have considered $G_2$ modular tensor categories. Kuperberg \cite{kuperberg:1996} gives a skein-theoretic description of quantum $G_2$ spider categories, which for generic $q$ agree with the fusion category of representations of the quantum group \cite{morrison/peters/snyder:2017}. In the case of $q$ a root of unity, Ostrik and Snyder have shown that the spider category agrees with the category of tilting modules for $G_2$ except for a few small roots of unity $q$ \cite{snyder/ostrik:2018talk}. As a consequence, the semisimplified versions of these categories obtained by quotienting by the ideal of negligible morphisms, again agree with the fusion category of representations of the quantum group.  These categories are in fact $C^{\ast}$-tensor categories and have been shown to have property (T) \cite{jones:2016}.

The paper is organised as follows. In Section \ref{sect:preliminaries} we present some preliminary material, including a discussion on the relation between spectral measures over certain different domains in Section \ref{sect:measures-different_domains} and a summary of relevant results for $G_2$ and its nimreps from \cite{evans/pugh:2012i}.

In Section \ref{sect:subgroupsG2} we discuss the finite subgroups of $G_2$, including their embeddings into the fundamental representations of $G_2$. We also give a general discussion of their spectral measures. In Section \ref{sect:results} we discuss each finite subgroup of $G_2$ individually, including determining all non-conjugate embeddings into the fundamental representations of $G_2$. For each such embedding we construct their McKay graphs, some of which have appeared before in \cite{he:2003} whilst the rest are new, and we determine their joint spectral measures.

\section{Preliminaries} \label{sect:preliminaries}
\subsection{Spectral measures over different domains} \label{sect:measures-different_domains}

The Weyl group of $G_2$ is the dihedral group $D_{12}$ of order 12.
As a subgroup of $GL(2,\bbZ)$, $D_{12}$ is generated by matrices $T_2$, $T_6$, of orders 2, 6 respectively, given by
\begin{equation} \label{T2,T6}
T_2 = \left( \begin{array}{cc} 0 & -1 \\ -1 & 0 \end{array} \right), \qquad T_6 = \left( \begin{array}{cc} 0 & 1 \\ -1 & 1 \end{array} \right),
\end{equation}
where the action of $D_{12}$ on $\bbT^2$ is given by $T(\omega_1,\omega_2) = (\omega_1^{a_{11}}\omega_2^{a_{12}},\omega_1^{a_{21}}\omega_2^{a_{22}})$, for $T = (a_{il}) \in D_{12}$. This action leaves $\chi_{\mu}(\omega_1,\omega_2)$ invariant, for any $\mu \in P_{++} = \{ (\mu_1,\mu_2) \in \mathbb{N}^2 | \, \mu_1 \geq \mu_2 \}$, the interior of the Weyl alcove for $G_2$.
Any $D_{12}$-invariant measure $\varepsilon_{\mu}$ on $\bbT^2$ yields a pushforward probability measure $\nu_{\mu}$ on $I_{\mu} = \chi_{\mu}( \bbT^2)\subset \bbR$ by
\begin{equation} \label{eqn:measures-T2-Ij_G2}
\int_{I_{\mu}} \psi(x) \mathrm{d}\nu_{\mu}(x) = \int_{\bbT^2} \psi(\chi_{\mu}(\omega_1,\omega_2)) \mathrm{d}\varepsilon_{\mu}(\omega_1,\omega_2),
\end{equation}
for any continuous function $\psi:I_{\mu} \rightarrow \mathbb{C}$, where $\mathrm{d}\varepsilon_{\mu}(\omega_1,\omega_2) = \mathrm{d}\varepsilon_{\mu}(g(\omega_1,\omega_2))$ for all $g \in D_{12}$.
There is a loss of dimension here, in the sense that the integral on the right hand side is over the two-dimensional torus $\bbT^2$, whereas on the right hand side it is over the interval $I_{\mu}$. Thus there is an infinite family of pullback measures $\varepsilon_{\mu}$ over $\bbT^2$ for any measure $\nu_{\mu}$ on $I_{\mu}$, that is, any $\varepsilon_{\mu}$ such that $\varepsilon_{\mu}(I_{\mu}^{-1}[x]) = \nu_{\mu}(x)$ for all $x \in I_{\mu}$ will yield the probability measure $\nu_{\mu}$ on $I_{\mu}$ as a pushforward measure by (\ref{eqn:measures-T2-Ij_G2}).
As in \cite{evans/pugh:2012i}, we instead work with an intermediate probability measure $\widetilde{\nu}_{\lambda,\mu}$ which lives over the joint spectrum $\mathfrak{D}_{\lambda,\mu} \subset I_{\lambda} \times I_{\mu} \subset \bbR^2$, for $\lambda,\mu \in P_{+}$, where there is no loss of dimension.

\begin{figure}[tb]
\begin{minipage}[t]{7.9cm}
\begin{center}
  \includegraphics[width=55mm]{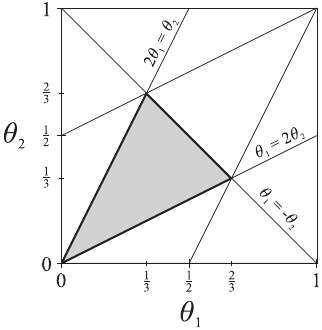}\\
 \caption{\small A fundamental domain of $\bbT^2/S_3$.} \label{fig:fund_domain-A2inT2}
\end{center}
\end{minipage}
\hfill
\begin{minipage}[t]{7.9cm}
\begin{center}
  \includegraphics[width=55mm]{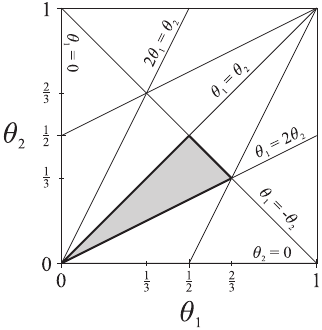}\\
 \caption{\small A fundamental domain $C$ of $\bbT^2/D_{12}$.} \label{fig:fund_domain-G2inT2}
\end{center}
\end{minipage}
\end{figure}

A fundamental domain $C$ of $\bbT^2$ under the action of the dihedral group $D_{12}$ is illustrated in Figure \ref{fig:fund_domain-G2inT2}, where the axes are labelled by the parameters $\theta_1$, $\theta_2$ in $(e^{2 \pi i \theta_1},e^{2 \pi i \theta_2}) \in \bbT^2$, which is a quotient of the fundamental domain of $\bbT^2/S_3$ illustrated in Figure \ref{fig:fund_domain-A2inT2} (see \cite{evans/pugh:2009v}) by the $\bbZ_2$-action given by the matrix -1.
Note that in Figure \ref{fig:fund_domain-G2inT2}, the lines $\theta_1=0$ and $\theta_2=0$ are also boundaries of copies of the fundamental domain $C$ under the action of $D_{12}$, whereas in Figure \ref{fig:fund_domain-A2inT2} they are not boundaries of copies of the fundamental domain under the action of $S_3$. The torus $\bbT^2$ contains 12 copies of $C$, so that
\begin{equation} \label{eqn:measureT2=12C}
\int_{\bbT^2} \phi(\omega_1,\omega_2) \mathrm{d}\varepsilon(\omega_1,\omega_2) = 12 \int_{C} \phi(\omega_1,\omega_2) \mathrm{d}\varepsilon(\omega_1,\omega_2),
\end{equation}
for any $D_{12}$-invariant function $\phi:\bbT^2 \rightarrow \mathbb{C}$ and $D_{12}$-invariant measure $\varepsilon$ over $\bbT^2$. The only fixed point of $\bbT^2$ under the action of $D_{12}$ is the point $(1,1)$.

Let $x_{\lambda} = \chi_{\lambda}(\omega_1,\omega_2)$ and let $\Psi_{\lambda,\mu}$ be the map $(\omega_1,\omega_2) \mapsto (x_{\lambda},x_{\mu})$. We denote by $\mathfrak{D}_{\lambda,\mu}$ the image of $\Psi_{\lambda,\mu}(C) (= \Psi_{\lambda,\mu}(\bbT^2))$ in $\bbR^2$.
Note that we can identify $\mathfrak{D}_{\lambda,\mu}$ with $\mathfrak{D}_{\mu,\lambda}$ by reflecting about the line $x_{\lambda} = x_{\mu}$.
Then the joint spectral measure $\widetilde{\nu}_{\lambda,\mu}$ is the measure on $\mathfrak{D}_{\lambda,\mu}$ uniquely determined by its cross-moments as in (\ref{eqn:cross_moments_sa_operators}).
Then there is a unique $D_{12}$-invariant pullback measure $\varepsilon$ on $\bbT^2$ such that
\begin{equation} \label{eqn:measures-T2-D_G2}
\int_{\mathfrak{D}_{\lambda,\mu}} \psi(x_{\lambda},x_{\mu}) \mathrm{d}\widetilde{\nu}_{\lambda,\mu}(x_{\lambda},x_{\mu}) = \int_{\bbT^2} \psi(\chi_{\lambda}(\omega_1,\omega_2),\chi_{\mu}(\omega_1,\omega_2)) \mathrm{d}\varepsilon_{\lambda,\mu}(\omega_1,\omega_2),
\end{equation}
for any continuous function $\psi:\mathfrak{D}_{\lambda,\mu} \rightarrow \mathbb{C}$.

Any probability measure on $\mathfrak{D}_{\lambda,\mu}$ yields a probability measure on the interval $I_{\lambda}$, given by the pushforward $(p_{\lambda})_{\ast}(\widetilde{\nu}_{\lambda,\mu})$ of the joint spectral measure $\widetilde{\nu}_{\lambda,\mu}$ under the orthogonal projection $p_{\lambda}$ onto the spectrum $\sigma(\lambda)$ (see \cite{evans/pugh:2012i} for more details).
Since the spectral measure $\nu_{\lambda}$ over $I_{\lambda}$ is also uniquely determined by its (one-dimensional) moments $\widetilde{\varsigma}_m = \int_{I_{\lambda}} x_{\lambda}^m \mathrm{d}\nu_{\lambda}(x_{\lambda})$ for all $m \in \mathbb{N}$, one could alternatively consider the moments in (\ref{eqn:cross_moments_sa_operators}) with $n=0$ to determine the measure $\nu_{\lambda}$ over $I_{\lambda}$.

Let
\begin{equation} \label{def:Dl}
C_k^W = \{ (e^{2 \pi i q_1/3(k+4)}, e^{2 \pi i q_2/3(k+4)}) \in \bbT^2 | \; q_1,q_2 = 0, 1, \ldots, 3k+11; \, q_1 + q_2 \equiv 0 \textrm{ mod } 3 \}
\end{equation}
which is the support (over $\bbT^2$) of the spectral measure of the nimrep graph $\mathcal{A}_k(G_2)$ associated to the trivial $G_2$ modular invariant at level $k$.
The following $G_2$-invariant measures will be useful later, c.f. \cite{evans/pugh:2010i}.
\begin{Def} \label{def:4measures}
Let $\omega = e^{2 \pi i/3}$, $\tau = e^{2 \pi i/n}$. We define the following measures on $\bbT^2$:
\begin{itemize}
\item[(1)] $\mathrm{d}_m \times \mathrm{d}_n$, where $\mathrm{d}_k$ is the uniform measure on the $k^{\mathrm{th}}$ roots of unity, for $k \in \mathbb{N}$.
\item[(2)] $\mathrm{d}^{(n)}$, the uniform measure on $C_n^W$ for $n \in \mathbb{N}$.
\item[(3)] $\mathrm{d}^{((n))}$, the uniform measure on the $S_3$-orbit of the points $(\tau, \tau)$, $(\overline{\omega} \, \overline{\tau}, \omega)$, $(\omega, \overline{\omega} \, \overline{\tau})$, for $n \in \mathbb{Q}$, $n \geq 2$.
\item[(4)] $\mathrm{d}^{(n,k)}$, the uniform measure on the $S_3$-orbit of the points $(\tau \, e^{2 \pi i k}, \tau)$, $(\tau, \tau \, e^{2 \pi i k})$, $(\overline{\omega} \, \overline{\tau}, \omega \, e^{2 \pi i k})$, $(\omega \, e^{2 \pi i k}, \overline{\omega} \, \overline{\tau})$, $(\overline{\omega} \, \overline{\tau} \, e^{-2 \pi i k}, \omega \, e^{-2 \pi i k})$, $(\omega \, e^{-2 \pi i k}, \overline{\omega} \, \overline{\tau} \, e^{-2 \pi i k})$, for $n,k \in \mathbb{Q}$, $n > 2$, $0 \leq k \leq 1/n$.
\end{itemize}
\end{Def}

\begin{figure}[tb]
\begin{minipage}[t]{7.5cm}
\begin{center}
  \includegraphics[width=55mm]{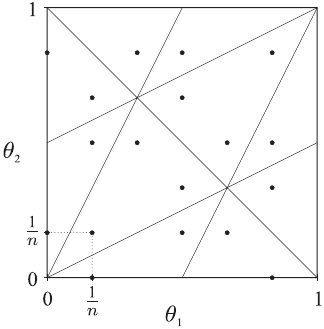}\\
 \caption{$\mathrm{Supp}(\mathrm{d}^{((n))})$} \label{fig:poly-15}
\end{center}
\end{minipage}
\hfill
\begin{minipage}[t]{7.5cm}
\begin{center}
  \includegraphics[width=55mm]{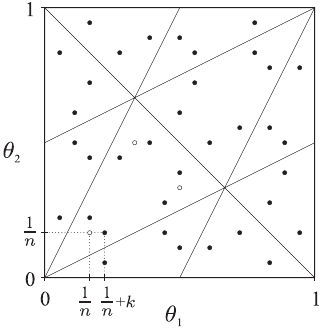}\\
 \caption{$\mathrm{Supp}(\mathrm{d}^{(n,k)})$} \label{fig:poly-16}
\end{center}
\end{minipage}
\end{figure}

The sets $\mathrm{Supp}(\mathrm{d}^{((n))})$, $\mathrm{Supp}(\mathrm{d}^{(n,k)})$ are illustrated in Figures \ref{fig:poly-15}, \ref{fig:poly-16} respectively, where $\mathrm{Supp}(\mathrm{d}\mu)$ denotes the set of points $(\theta_1,\theta_2) \in [0,1]^2$ such that $(e^{2 \pi i \theta_1}, e^{2 \pi i \theta_2})$ is in the support of the measure $\mathrm{d}\mu$. The white circles in Figure \ref{fig:poly-16} denote the points given by the measure $\mathrm{d}^{((n))}$. The cardinality $|\mathrm{Supp}(\mathrm{d}_m \times \mathrm{d}_n)|$ of $\mathrm{Supp}(\mathrm{d}_m \times \mathrm{d}_n)$ is $mn$, whilst $|\mathrm{Supp}(\mathrm{d}^{(n)})| = |D_n| = 3n^2$ was shown in \cite[Section 7.1]{evans/pugh:2009v}. For $n > 2$ and $0 < k < 1/n$, $|\mathrm{Supp}(\mathrm{d}^{((n))})| = 18$, whilst $|\mathrm{Supp}(\mathrm{d}^{(n,k)})| = 36$. The cardinalities of the other sets are $|\mathrm{Supp}(\mathrm{d}^{(n,0)})| = |\mathrm{Supp}(\mathrm{d}^{(n,1/n)})| = 18$ for $n > 2$, and $|\mathrm{Supp}(\mathrm{d}^{((2))})| = 9$.
Some relations between these measures are given in \cite[Section 2]{evans/pugh:2010i}.

\subsection{Spectral measures for $G_2$} \label{sect:spec_measure-G2}

Here we review the results, determined in \cite{evans/pugh:2012i}, for the spectral measures for $G_2$.
Let $\rho_1$, $\rho_2$ denote the fundamental representations of $G_2$ of dimensions 7, 14 respectively.
The restrictions of the characters $\chi_{\rho_j}$ of the fundamental representations of $G_2$ to $\bbT^2$ yield maps from the torus to the interval $I_j = \chi_{\rho_j}(\bbT^2) \subset \bbR$:
\begin{align*}
\chi_{\rho_1}(\omega_1,\omega_2) & = 1 + \omega_1 + \omega_1^{-1} + \omega_2 + \omega_2^{-1} + \omega_1\omega_2^{-1} + \omega_1^{-1}\omega_2 \\
& = 1 + 2\cos(2\pi\theta_1) + 2\cos(2\pi\theta_2) + 2\cos(2\pi(\theta_1-\theta_2)), \\
\chi_{\rho_2}(\omega_1,\omega_2) & = \chi_{\rho_1}(\omega_1,\omega_2) + 1 + \omega_1\omega_2 + \omega_1^{-1}\omega_2^{-1} + \omega_1^2\omega_2^{-1} + \omega_1^{-2}\omega_2 + \omega_1\omega_2^{-2} + \omega_1^{-1}\omega_2^2 \\
= \chi_{\rho_1}&(\omega_1,\omega_2) + 1 + 2\cos(2\pi(\theta_1+\theta_2)) + 2\cos(2\pi(2\theta_1-\theta_2)) + 2\cos(2\pi(\theta_1-2\theta_2)),
\end{align*}
where $\omega_j = e^{2\pi i \theta_j} \in \bbT$ for $\theta_j \in [0,1]$, $j=1,2$.

\begin{figure}[tb]
\begin{center}
  \includegraphics[width=60mm]{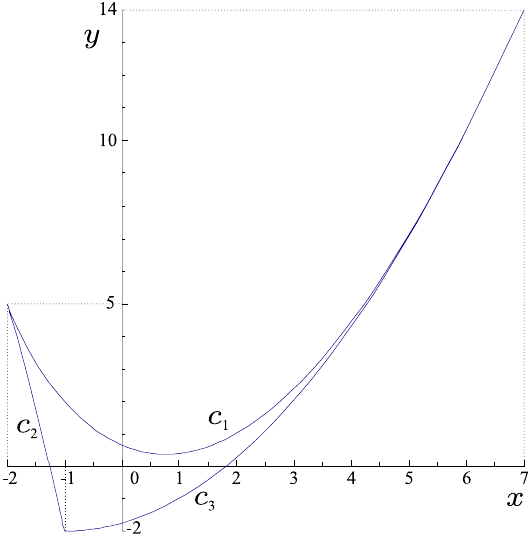}\\
 \caption{The domain $\mathfrak{D} = \Psi(C)$.} \label{fig:DomainD-G2}
\end{center}
\end{figure}

Let
\begin{equation} \label{eqn:x,y-G2}
x := \chi_{\rho_1}(\omega_1,\omega_2), \qquad y := \chi_{\rho_2}(\omega_1,\omega_2),
\end{equation}
and denote by $\Psi$ the map $\Psi_{(1,0),(1,1)}: (\omega_1,\omega_2) \mapsto (x,y)$.
The image $\mathfrak{D} = \Psi(C)$ of the fundamental domain $C$ of $\bbT^2/D_{12}$ (illustrated in Figure \ref{fig:fund_domain-G2inT2}) under $\Psi$ is illustrated in Figure \ref{fig:DomainD-G2}, where the boundaries of $C$ given by $\theta_1 = 2\theta_2$, $\theta_1 = -\theta_2$, $\theta_1 = \theta_2$ yield the curves $c_1$, $c_2$, $c_3$ respectively. These curves are given by \cite{uhlmann/meinel/wipf:2007}
\begin{align*}
c_1: && y & = -5(x+1)+2(x+2)^{3/2}, \qquad x \in [-2,7], \\
c_2: && y & = -5(x+1)-2(x+2)^{3/2}, \qquad x \in [-2,-1], \\
c_3: && 4y & = x^2+2x-7, \hspace{28mm} x \in [-1,7].
\end{align*}
The fixed point $(1,1)$ of $\bbT^2$ under the action of $D_{12}$ maps to 7, 14 in the intervals $I_1$, $I_2$ respectively.

Under the change of variables (\ref{eqn:x,y-G2}) the Jacobian is given by \cite{evans/pugh:2012i}
\begin{equation} \label{eqn:J[theta]-D12}
\begin{split}
J & = 8 \pi^2 (\cos(2 \pi (2\theta_1 + \theta_2)) + \cos(2 \pi (\theta_1 - 3\theta_2)) + \cos(2 \pi (3\theta_1 - 2\theta_2)) \\
& \qquad - \cos(2 \pi (\theta_1 + 2\theta_2)) - \cos(2 \pi (3\theta_1 - \theta_2)) - \cos(2 \pi (2\theta_1 - 3\theta_2))).
\end{split}
\end{equation}
The Jacobian is real and vanishes in $\bbT^2$ only on the boundaries of the images of the fundamental domain $C$ under $D_{12}$.
Again, $J^2$ can be written in terms of the $D_{12}$-invariant elements $x$, $y$ as $J^2 = (4x^3-x^2-2x-10xy-y^2-10y+7)(x^2+2x-7-4y)$ (see also \cite{uhlmann/meinel/wipf:2007}), which is non-negative since $J$ is real. We write $J$ in terms of $x$ and $y$ as
\begin{eqnarray} \label{eqn:J[x,y]-D12}
|J| & = & 4 \pi^2 \sqrt{(4x^3-x^2-2x-10xy-y^2-10y+7)(x^2+2x-7-4y)}.
\end{eqnarray}

The joint spectral measure over $\mathfrak{D}$ is then $|J(x,y)| \mathrm{d}x \, \mathrm{d}y/16\pi^4$, which is the reduced Haar measure on $G_2$ \cite[Sect. 6.3]{uhlmann/meinel/wipf:2007}.

\subsection{Spectral measures for nimrep graphs associated to $G_2$ modular invariants} \label{sect:spec_measure-nimrepsG2}

Suppose $G$ is the nimrep associated to a $G_2$ braided subfactor at some finite level $k$ with vertex set $G_0$.
We define a state $\varphi$ on $\ell^2(G_0)$ by $\varphi( \, \cdot \, ) = \langle \,\cdot \, e_{\ast}, e_{\ast} \rangle$, where $e_{\ast}$ is the basis vector in $\ell^2(G_0)$ corresponding to the distinguished vertex $\ast$ with lowest Perron-Frobenius weight.

If we consider the nimrep graphs $G_{\lambda}$, $G_{\mu}$, which have joint spectrum $\mathfrak{D}_{\lambda,\mu}$, then the $m,n^{\mathrm{th}}$ cross moment $\varsigma_{m,n} = \varphi(G_{\lambda}^m G_{\mu}^n) = \int_{\mathfrak{D}_{\lambda,\mu}} x^m y^n \mathrm{d}\widetilde{\nu}(x,y)$, where $x=x_{\lambda}$, $y=x_{\mu}$, is given by $\langle G_{\lambda}^m G_{\mu}^n e_{\ast}, e_{\ast} \rangle$.
Let $\beta_{\lambda}^{(\nu)} = \chi_{\lambda}(t_{\nu})$ be the eigenvalues of $G_{\lambda}$, where $t_{\nu}=(\exp(\xi(\nu_1+1)),\exp(-3\xi (\nu_2+1)))$ for $\xi = 6\pi i/(k+4)$, with corresponding eigenvectors $\psi^{(\nu)}$ (note that the eigenvectors of $G_{\lambda}$ are the same for all $\lambda$, since the $G_{\lambda}$ are a family of commuting matrices). Each eigenvalue $\beta_{\lambda}^{(\mu)}$ is also given by a ratio of the modular $S$-matrix for $G_2$ at level $k$, $\beta_{\lambda}^{(\mu)} = S_{\lambda\mu}/S_{0\mu}$, where $\mu \in \mathrm{Exp}(G) \subset P^k_{+} = \{ (\lambda_1,\lambda_2) | \, \lambda_1,\lambda_2 \geq 0; \lambda_1 + 2\lambda_2 \leq k \}$ are given by the modular invariant $Z$.
Then $G_{\lambda}^m G_{\mu}^n = \mathcal{U} \Lambda_{\lambda}^m \Lambda_{\mu}^n \mathcal{U}^{\ast}$, where $\Lambda_{\lambda}$ is the diagonal matrix with the eigenvalues $\beta_{\lambda}^{(\nu)}$ on the diagonal, and $\mathcal{U}$ is the matrix whose columns are given by the eigenvectors $\psi^{(\nu)}$, so that
\begin{equation}\label{eqn:moments-nimrep-G2}
\varsigma_{m,n} \;\; = \;\; \langle \mathcal{U} \Lambda_{\lambda}^m \Lambda_{\mu}^n \mathcal{U}^{\ast} e_{\ast}, e_{\ast} \rangle \;\; = \;\; \langle \Lambda_{\lambda}^m \Lambda_{\mu}^n \mathcal{U}^{\ast} e_{\ast}, \mathcal{U}^{\ast} e_{\ast} \rangle \;\; = \;\; \sum_{\nu} (\beta_{\lambda}^{(\nu)})^m (\beta_{\mu}^{(\nu)})^n |\psi^{(\nu)}_{\ast}|^2,
\end{equation}
where $\psi^{(\nu)}_{\ast} = \mathcal{U}^{\ast} e_{\ast}$ is the entry of the eigenvector $\psi^{(\nu)}$ corresponding to the distinguished vertex $\ast$.
Then there is a $D_{12}$-invariant measure $\varepsilon$ over $\bbT^2$ such that
$$\varsigma_{m,n} = \int_{\bbT^2} \chi_{\lambda}(\omega_1,\omega_2)^m \chi_{\mu}(\omega_1,\omega_2)^n \mathrm{d}\varepsilon(\omega_1,\omega_2),$$
for all $\lambda$, $\mu$.

Note from (\ref{eqn:moments-nimrep-G2}) that the measure $\varepsilon$ is a discrete measure which has weight $|\psi^{(\nu)}_{\ast}|^2$ at the points $g(t_{\nu}) \in \bbT^2$ for $g \in D_{12}$, $\nu \in \mathrm{Exp}(G)$, and zero everywhere else. Thus the measure $\varepsilon$ does not depend on the choice of $\lambda$, $\mu$, so that the spectral measure over $\bbT^2$ is the same for any pair $(G_{\lambda},G_{\mu})$, even though the corresponding measures over $\mathfrak{D}_{\lambda,\mu} \subset \bbR^2$, and indeed the subsets $\mathfrak{D}_{\lambda,\mu}$ themselves, are different for each such pair. The same result holds for the spectral measure over $\bbT^2$ of a finite subgroup of $G_2$.

\section{The finite subgroups of $G_2$} \label{sect:subgroupsG2}

The classification of finite subgroups of $G_2$ is due to \cite{wales:1970, cohen/wales:1983} (see also \cite{greiss:1995, he:2003}).

The reducible (i.e. block-diagonalizable) finite subgroups of $G_2$ are the finite discrete subgroups of $SU(2) \times SU(2)$ and $SU(3)$ \cite{wales:1970}. These subgroups are thus well known, and the corresponding spectral measures can be obtained from \cite{banica/bisch:2007, evans/pugh:2009v, evans/pugh:2010i}.

The irreducible finite subgroups of $G_2$, of which there are seven up to conjugacy in $G_2$ (or equivalently, up to conjugacy in $GL(V)$, where $V$ is the natural 7-dimensional module for $O(7,\mathbb{C})$ \cite[Corollary 1]{greiss:1995}), can be further classified into two types, primitive and imprimitive, where a linear group $\Gamma \subset GL(V)$ is imprimitive if there is a non-trivial decomposition $V=\bigoplus_i V_i$ such that $\Gamma$ permutes the $V_i$. There are two imprimitive finite subgroups and five primitive ones.
These finite subgroups are listed in Table \ref{Table:subgroupsG2}, where type denotes whether an irreducible subgroup is primitive (P) or imprimitive (I).

\renewcommand{\arraystretch}{1}

\begin{table}[tb]
\begin{center}
\begin{tabular}{|c|c|c|} \hline
Subgroup $\Gamma \subset G_2$ & Type & $|\Gamma|$ \\
\hline\hline finite subgroups of $SU(2) \times SU(2)$, $SU(3)$ & - & - \\
\hline $PSL(2;7) \cong GL(3;2) \cong \Sigma (168) \subset SU(3)$ & I & 168 \\
\hline $PSL(2;7) \rtimes \bbZ_2^3$ & I & 1344 \\
\hline $PGL(2;7)$ & P & 336 \\
\hline $PSL(2;8)$ & P & 504 \\
\hline $PSL(2;13)$ & P & 1092 \\
\hline $PU(3;3) \cong G_2(2)'$ & P & 6048 \\
\hline $G_2(2)$ & P & 12096 \\
\hline
\end{tabular} \\
\caption{Finite subgroups of $G_2$.} \label{Table:subgroupsG2}
\end{center}
\end{table}

The McKay graph $\mathcal{G}^{\rho}_{\Gamma}$ is the the fusion graph of the irreducible representation $\rho$ of $\Gamma$ acting on the irreducible representations of $\Gamma$.
This graph determines the Bratteli diagram for the tower of relative commutants of the subfactor $P^{\Gamma} \subset (M_n \otimes P)^{\Gamma}$, where $n$ is the dimension the representation $\rho$ and $P$ is the type $\mathrm{II}_1$ factor $\bigotimes_{n=1}^{\infty} M_n$ \cite[$\S$VI]{wassermann:1988}. This graph is not however the principal graph of this subfactor as it is not bipartite. The prinicpal graph is rather an unfolded version of the McKay graph $\mathcal{G}^{\rho}_{\Gamma}$, with adjacency matrix given by $\left( \begin{array}{cc} 0 & \Delta^{\rho}_{\Gamma} \\ \Delta^{\rho}_{\Gamma} & 0 \end{array} \right)$, where $\Delta^{\rho}_{\Gamma}$ is the adjacency matrix of the (folded) graph $\mathcal{G}^{\rho}_{\Gamma}$.

We will consider the (joint) spectral measure for the McKay graphs $\mathcal{G}^j_{\Gamma} := \mathcal{G}^{\varrho_j}_{\Gamma}$ associated to a finite subgroup $\Gamma \subset G_2$, where $\varrho_j$ are the restrictions of the fundamental representations $\rho_j$ of $G_2$ to $\Gamma$, $j=1,2$.
We will consider all possible embeddings of the subgroup in $G_2$. Any two such embeddings are conjguate in $G_2$ if and only if they afford the same character on the seven-dimensional representation $\rho_1$ \cite[Corollary 1]{greiss:1995}. In some cases there is more than one non-conjugate embedding of the subgroup in $G_2$.
In these cases the restricted representation $\varrho_1$ is not necessarily irreducible.
In all cases, even for irreducible $\varrho_1$, the restriction $\varrho_2$ of the 14-dimensional representation is not necessarily irreducible.

We use the following methods to determine embeddings $\varrho_1$ of $\Gamma$ in $G_2$.
First, take a seven-dimensional (not necessarily irreducible) representation $\gamma_1$ of $\Gamma$.
The Kronecker square of $\rho_1$ decomposes into irreducible representations of $G_2$ as $\rho_1^2 = \mathrm{id}_{G_2} + \rho_1 + \rho_2 + \lambda_{(2,0)}$, where $\lambda_{(2,0)}$ has dimension 27. The Kronecker square of $\varrho_1$ is obtained by restricting this decomposition to $\Gamma$, and we see that $\varrho_1$ appears in the decomposition of $\varrho_1^2$ into irreducible representations of $\Gamma$. Thus, if $\gamma_1$ is not contained in the decomposition of $\gamma_1^2$ into irreducible representations of $\Gamma$, we can eliminate $\gamma_1$ as a possible restriction $\varrho_1$ of $\rho_1$. The decomposition of $\gamma_1^2$ into irreducible representations can be obtained using the character table for $\Gamma$, by decomposing the character $\chi_{\gamma_1^2} = \chi_{\gamma_1}^2$ of $\gamma_1^2$ into the characters of the irreducible representations $\lambda$ of $\Gamma$: $\chi_{\gamma_1}^2 = \sum_{\lambda} a_{\lambda} \chi_{\lambda}$, where $a_{\lambda} = \langle \gamma_1^2, \lambda \rangle/|\Gamma| = \sum_{g \in \Gamma} \chi_{\gamma_1^2}(g)\chi_{\lambda}(g)/|\Gamma|$.

We next consider the eigenvalues of the representation matrices of $\Gamma$. If the elements in a conjugacy class $C_n$ of $\Gamma$ have order $n$, then $(C_n)^n = Z(\Gamma)$, where $Z(\Gamma)$ is the center of $\Gamma$. If the center is trivial, which is the case for all subgroups of $G_2$, then the eigenvalues $\xi$ of the matrices representing these elements must satisfy $\xi^n = 1$. Since $\chi_{\lambda}(\Gamma_j)$ is the sum of the eigenvalues $\xi$, it is usually possible to write down the complete set of eigenvalues from the information provided by the character table of $\Gamma$ and the fact that the eigenvalues must be powers of $n^{\mathrm{th}}$ roots of unity.
Where there is some ambiguity, we can pin down the correct choice for the set of eigenvalues from the following considerations. Suppose there is ambiguity regarding the eigenvalues of the conjugacy class $C_{mn}$ whose elements have order $mn$, $m,n \in \mathbb{N}$. If there is only one conjugacy class $C_n$ whose elements have order $n$, then for $g \in C_{mn}$, $g^m \in C_n$, and since $g, g^m$ commute, their corresponding representation matrices can be simultaneously diagonalised, and thus the eigenvalues of $C_n$ must be $m^{\mathrm{th}}$ powers of those for $C_{mn}$. Suppose now that there is more than one conjugacy class $C_n^{(j)}$ whose elements have order $n$. Since $g^m$ are all conjugate for conjugate $g$, we see that there exists a $j$ such that $g^m \in C_n^{(j)}$ for all $g \in C_{mn}$. It turns out in all the cases considered here that there is only one consistent choice of $j$ such that the eigenvalues of $C_n^{(j)}$ are $m^{\mathrm{th}}$ powers of those for $C_{mn}$ for all (irreducible) representations.

As was shown in \cite[$\S$4]{evans/pugh:2009v}, the eigenvalues of the representation matrices of $\Gamma$ can be written in the form $\chi_{\varrho_j}(C) = \mathrm{Tr}(\varrho_j(g))$, where $g$ is any element of the conjugacy class $C$ of $\Gamma$.
Every element $g \in \Gamma$ is conjugate to an element $d$ in the maximal torus of $G_2$, i.e. $\varrho_j(h^{-1}gh) = \varrho_j(d) = (\varrho_j|_{\bbT^2})(t_1,t_2)$ for some $(t_1,t_2) \in \bbT^2$, for $j=1,2$, where $\varrho_j|_{\bbT^2}$ is given by \cite{evans/pugh:2012i}
\begin{align}
(\rho_1|_{\bbT^2})&(t_1,t_2) = \textrm{diag}(D(t_1), D(t_2^{-1}), D(t_1^{-1}t_2), 1), \label{eqn:restrict_rho1G2_to_T2} \\
(\rho_2|_{\bbT^2})&(t_1,t_2) \nonumber \\
&= \textrm{diag}(D(t_1), D(t_2^{-1}), D(t_1^{-1}t_2), D(1), D(t_1t_2), D(t_1^2t_2^{-1}), D(t_1^{-1}t_2^2)), \label{eqn:restrict_rho2G2_to_T2}
\end{align}
where $D(t_i) = \left( \begin{array}{cc} \mathrm{Re}(t_i) & -\mathrm{Im}(t_i) \\ \mathrm{Im}(t_i) & \mathrm{Re}(t_i) \end{array} \right)$ for $t_i \in \bbT$.
Now $\mathrm{Tr}(\varrho_j(g)) = Tr(\varrho_j(d)) = \Phi_j(t_1,t_2)$, thus the eigenvalues of $\varrho_j(g)$ are all of the form $\Phi_j(\omega_1,\omega_2)$ for $\omega_1,\omega_2 \in \bbT$, and hence its spectrum is contained in the interval $I_j$.
As shown in \cite[Sections 3,4]{evans/pugh:2012i} the spectrum of the fundamental representation $\rho_j$ of $G_2$, and its restriction to $\bbT^2$, is the whole of the interval $I_j = \chi_{\rho_j}(\bbT^2)$, for $j=1,2$.
Thus the support of the spectral measure $\mu_{\Delta_j}$ of $\Delta_j = \Delta_{\mathcal{G}^j_{\Gamma}}$, the adjacency matrix of $\mathcal{G}^j_{\Gamma}$, is contained in $I_j$ when $\Gamma$ is $G_2$ or one of its finite subgroups.
Then for $\Gamma \subset G_2$, the eigenvalues of every group element in $\varrho_1$ are necessarily of the form $\mathcal{E}_{t_1,t_2} := \{ 1,t_1,t_1^{-1},t_2,t_2^{-1},t_1t_2^{-1},t_1^{-1}t_2 \}$, where $t_i \in \bbT$. Thus, by \cite[Proposition]{king/toumazet/wybourne:1999}, if the eigenvalues of the group elements in the representation $\gamma_1$ have this form then $\gamma_1$ is a restriction $\varrho_1$ of the seven-dimensional fundamental representation $\rho_1$ of $G_2$ to $\Gamma$.

We now turn to consider the possible restrictions $\varrho_2$ of the fourteen-dimensional fundamental representation $\rho_2$ of $G_2$ to $\Gamma$, for fixed $\varrho_1$.
By dimension considerations one can determine the possible candidates for $\varrho_2$ from the Kronecker square of $\varrho_1$.
Let $\gamma_2$ be such a candidate, i.e. a fourteen-dimensional representation such that $\varrho_1^2 = \mathrm{id}_{\Gamma} + \varrho_1 + \gamma_2 + \lambda$, where $\lambda$ is (necessarily) some 27-dimensional representation of $\Gamma$.
We make a choice of pair $(t_1^{C},t_2^{C})$ from the set $X_C$ of eigenvalues of group elements (from the conjugacy class $C$) in $\varrho_1$ such that $\mathcal{E}_{t_1^C,t_2^C} = X_C$.
Note that the choice of $(t_1^C,t_2^C)$ such that $\mathcal{E}_{t_1^C,t_2^C} = X_C$ is not unique. However, any other pair $(\tilde{t}_1^C,\tilde{t}_2^C)$ such that $\mathcal{E}_{\tilde{t}_1^C,\tilde{t}_2^C} = X_C$ will appear in the orbit of $(t_1^C,t_2^C)$ under the action of the Weyl group $D_{12}$ of $G_2$, where the action of $D_{12}$ on $\bbT^2$ is given in Section \ref{sect:measures-different_domains}. Thus $\Phi_j(\tilde{t}_1^C,\tilde{t}_2^C) = \Phi_j(t_1^C,t_2^C)$ for $j=1,2$.
Then one checks that $\Phi_2(t_1^C,t_2^C) = \chi_{\gamma_2}(C)$ for each conjugacy class $C$, in which case $\gamma_2$ is indeed a restriction $\varrho_2$ of the fourteen-dimensional fundamental representation $\rho_2$ of $G_2$ to $\Gamma$.

The following result present the classification of all non-conjugate embeddings of the finite subgroups of $G_2$ in $G_2$. For the details of the embeddings for each individual subgroup see the corresponding subsections of Section \ref{sect:results}.

\begin{Thm}
Up to conjugacy in $G_2$ the number of faithful embeddings of each finite subgroup of $G_2$ in the seven-dimensional fundamental representation of $G_2$ is:
\begin{itemize}
\item For $PU(3;3)$ and $G_2(2)$ there is only one embedding;
\item For $PSL(2;7)$, $\bbZ_2^3 \cdot PSL(2;7)$, $PGL(2;7)$ and $PSL(2;13)$ there are two non-conjugate embeddings;
\item For $PSL(2;8)$ there are three non-conjugate embeddings.
\end{itemize}
\end{Thm}

\subsection{Spectral measures for finite subgroups of $G_2$} \label{sect:spec_measure-subgroupsG2}

For finite groups, the analogue of the $S$-matrix of Section \ref{sect:spec_measure-nimrepsG2}, which simultaneously diagonalizes the representations of $\Gamma$, is essentially the character table of $\Gamma$. This $S$-matrix, with rows, columns labelled by the irreducible characters, conjugacy classes respectively of $\Gamma$, is given by $S_{\lambda,C} = \sqrt{|C|} \chi_{\lambda}(C)/\sqrt{|\Gamma|}$, for representations $\lambda$ and conjugacy classes $C$. Then the analogue of \eqref{eqn:moments-nimrep-G2} for finite groups is that the $m,n^{\mathrm{th}}$ moment $\varsigma_{m,n}$ is given over $\mathfrak{D}$ by \cite{kawai:1989} (c.f. \cite[Section 4]{evans/pugh:2009v} for the case of finite subgroups of $SU(2)$)
\begin{equation} \label{eqn:moments-subgroupG2}
\varsigma_{m,n} \; = \; \int_{\mathfrak{D}} x^m y^n \mathrm{d}\nu(x,y) \; = \; \sum_C \frac{|C|}{|\Gamma|} \chi_{\varrho_1} (C)^m \chi_{\varrho_2} (C)^n.
\end{equation}
There is an analogous statement to (\ref{eqn:moments-subgroupG2}) for the joint spectral measure $\nu_{\lambda,\mu}$ over $\mathfrak{D}_{\lambda,\mu}$ for any irreducible representations $\lambda$, $\mu$ of $\Gamma$. The weight on the right hand side will again be $|C|/|\Gamma|$, since the same $S$-matrix simultaneously diagonalises all the representations of $\Gamma$. Since this weight does not depend on the representations $\lambda$, $\mu$, we see that the $D_{12}$-invariant pullback measure $\varepsilon$ over $\bbT^2$ will be the same for any joint spectral measure $\nu_{\lambda,\mu}$. This is analogous to the situation for nimrep graphs discussed in Section \ref{sect:spec_measure-nimrepsG2}.

We wish to compute `inverse' maps $\widetilde{\Psi}: \mathfrak{D} \rightarrow \bbT^2$ such that $\Psi \circ \widetilde{\Psi} = \mathrm{id}$.
The following equation can easily be checked by substituting in $x=\Phi_1(\omega_1,\omega_2)$, $y=\Phi_2(\omega_1,\omega_2)$:
$$(\omega_j+\omega_j^{-1})^3 + (1-x)(\omega_j+\omega_j^{-1})^2 + (y-2)(\omega_j+\omega_j^{-1}) +2y-x^2+2x-1 = 0,$$
where $j=1,2$.
Solving this cubic in $\omega_j+\omega_j^{-1} = 2\cos(\vartheta_j)$, we obtain for $l \in \{ 0,1,2 \}$
$$\vartheta_j(l) = \cos^{-1}\left( \frac{1}{6} \left( x-1 + 2^{-1/3} \epsilon_l P + 2^{1/3} \overline{\epsilon_l} (x^2-2x+7-3y)P^{-1} \right) \right)$$
where $2^{1/3}$ takes a real value, $\epsilon_l = e^{2 \pi i l/3}$ and $P = (2x^3+21x^2-30x+7-45y-9xy + \sqrt{(4x^3-x^2-2x-10xy-y^2-10y+7)(x^2+2x-7-4y)} \, )^{1/3}$. We note that for the roots of a cubic equation it does not matter whether the square root in $P$ is taken to be positive or negative.
Then we set $\Psi_{l,l'}(x,y) = (e^{\vartheta_1(l)i},e^{\vartheta_2(l')i})$, and we have that $\Psi(\Psi_{l,l'}(x,y)) = (x,y)$ for some $l,l' \in \{ 0,1,2 \}$. The particular choice of pair $l,l'$ such that the equality $\Psi \circ \Psi_{l,l'} = \mathrm{id}$ is satisfied depends on $x,y$, but it is easy to check (eg. using Mathematica) whether a given choice satisfies this equality for any of the examples we consider.
We present in Table \ref{Table:subgroupsG2-orbits(theta1,theta2)} the values of the eigenvalues $(\chi_{\varrho_1}(C), \chi_{\varrho_2}(C)) = (x,y) \in \mathfrak{D}$ which will appear for the finite subgroups of $G_2$, and (the orbits under $D_{12}$ of) the corresponding points $(\theta_1,\theta_2) \in [0,1]^2$ such that $\Psi(e^{2\pi i \theta_1},e^{2\pi i \theta_2}) = (x,y)$.

\renewcommand{\arraystretch}{1.4}

\begin{table}[tbp]
\begin{center}
\begin{tabular}{|c|c|c|} \hline
$(x,y) \in \mathfrak{D}$ & Orbit of $(\theta_1,\theta_2) \in [0,1]^2$ & $\displaystyle \frac{J^2(\theta_1,\theta_2)}{16\pi^4}$ \\
\hline $(7,14)$ & $(0,0)$ & 0 \\
\hline $(-2,5)$ & $\left(\frac{1}{3},\frac{2}{3}\right), \left(\frac{2}{3},\frac{1}{3}\right)$ & 0 \\
\hline $(-1,-2)$ & $\left(0,\frac{1}{2}\right), \left(\frac{1}{2},\frac{1}{2}\right), \left(\frac{1}{2},0\right)$ & 0 \\
\hline $(1,-1)$ & $\left(0,\frac{1}{3}\right), \left(\frac{1}{3},\frac{1}{3}\right), \left(\frac{1}{3},0\right), \left(0,\frac{2}{3}\right), \left(\frac{2}{3},\frac{2}{3}\right), \left(\frac{2}{3},0\right)$ & 0 \\
\hline $(3,2)$ & $\left(0,\frac{1}{4}\right), \left(\frac{1}{4},\frac{1}{4}\right), \left(\frac{1}{4},0\right), \left(0,\frac{3}{4}\right), \left(\frac{3}{4},\frac{3}{4}\right), \left(\frac{3}{4},0\right)$ & 0 \\
\hline $(-1,2)$ & $\left(\frac{1}{4},\frac{1}{2}\right), \left(\frac{1}{2},\frac{1}{4}\right), \left(\frac{1}{4},\frac{3}{4}\right), \left(\frac{3}{4},\frac{1}{2}\right), \left(\frac{1}{2},\frac{3}{4}\right), \left(\frac{3}{4},\frac{1}{4}\right)$ & 0 \\
\hline $(2,1)$ & $\left(\frac{1}{6},\frac{1}{3}\right), \left(\frac{1}{3},\frac{1}{6}\right), \left(\frac{1}{6},\frac{5}{6}\right), \left(\frac{5}{6},\frac{2}{3}\right), \left(\frac{2}{3},\frac{5}{6}\right), \left(\frac{5}{6},\frac{1}{6}\right)$ & 0 \\
\hline $(-1,1)$ & $\left(\frac{1}{6},\frac{1}{2}\right), \left(\frac{1}{2},\frac{1}{3}\right), \left(\frac{1}{3},\frac{5}{6}\right), \left(\frac{5}{6},\frac{1}{2}\right), \left(\frac{1}{2},\frac{2}{3}\right), \left(\frac{2}{3},\frac{1}{6}\right),$ & 36 \\
& $\left(\frac{1}{2},\frac{1}{6}\right), \left(\frac{1}{3},\frac{1}{2}\right), \left(\frac{5}{6},\frac{1}{3}\right), \left(\frac{1}{2},\frac{5}{6}\right), \left(\frac{2}{3},\frac{1}{2}\right), \left(\frac{1}{6},\frac{2}{3}\right)$ & \\
\hline $(0,0)$ & $\left(\frac{1}{7},\frac{3}{7}\right), \left(\frac{3}{7},\frac{2}{7}\right), \left(\frac{2}{7},\frac{6}{7}\right), \left(\frac{6}{7},\frac{4}{7}\right), \left(\frac{4}{7},\frac{5}{7}\right), \left(\frac{5}{7},\frac{1}{7}\right),$ & $49$ \\
& $\left(\frac{3}{7},\frac{1}{7}\right), \left(\frac{2}{7},\frac{3}{7}\right), \left(\frac{6}{7},\frac{2}{7}\right), \left(\frac{4}{7},\frac{6}{7}\right), \left(\frac{5}{7},\frac{4}{7}\right), \left(\frac{1}{7},\frac{5}{7}\right)$ & \\
\hline $(1,0)$ & $\left(\frac{1}{8},\frac{3}{8}\right), \left(\frac{3}{8},\frac{1}{4}\right), \left(\frac{1}{4},\frac{7}{8}\right), \left(\frac{7}{8},\frac{5}{8}\right), \left(\frac{5}{8},\frac{3}{4}\right), \left(\frac{3}{4},\frac{1}{8}\right),$ & 32 \\
& $\left(\frac{3}{8},\frac{1}{8}\right), \left(\frac{1}{4},\frac{3}{8}\right), \left(\frac{7}{8},\frac{1}{4}\right), \left(\frac{5}{8},\frac{7}{8}\right), \left(\frac{3}{4},\frac{5}{8}\right), \left(\frac{1}{8},\frac{3}{4}\right)$ & \\
\hline $(-1,0)$ & $\left(\frac{1}{8},\frac{1}{2}\right), \left(\frac{1}{2},\frac{3}{8}\right), \left(\frac{3}{8},\frac{7}{8}\right), \left(\frac{7}{8},\frac{1}{2}\right), \left(\frac{1}{2},\frac{5}{8}\right), \left(\frac{5}{8},\frac{1}{8}\right),$ & 32 \\
& $\left(\frac{1}{2},\frac{1}{8}\right), \left(\frac{3}{8},\frac{1}{2}\right), \left(\frac{7}{8},\frac{3}{8}\right), \left(\frac{1}{2},\frac{7}{8}\right), \left(\frac{5}{8},\frac{1}{2}\right), \left(\frac{1}{8},\frac{5}{8}\right)$ & \\
\hline $(-p,p+q+1)$ & $\left(\frac{1}{9},\frac{4}{9}\right), \left(\frac{4}{9},\frac{1}{3}\right), \left(\frac{1}{3},\frac{8}{9}\right), \left(\frac{8}{9},\frac{5}{9}\right), \left(\frac{5}{9},\frac{2}{3}\right), \left(\frac{2}{3},\frac{1}{9}\right),$ & $9a_1$ \\
& $\left(\frac{4}{9},\frac{1}{9}\right), \left(\frac{1}{3},\frac{4}{9}\right), \left(\frac{8}{9},\frac{1}{3}\right), \left(\frac{5}{9},\frac{8}{9}\right), \left(\frac{2}{3},\frac{5}{9}\right), \left(\frac{1}{9},\frac{2}{3}\right)$ & \\
\hline $(-q,1-p)$ & $\left(\frac{1}{9},\frac{1}{3}\right), \left(\frac{1}{3},\frac{2}{9}\right), \left(\frac{2}{9},\frac{8}{9}\right), \left(\frac{8}{9},\frac{2}{3}\right), \left(\frac{2}{3},\frac{7}{9}\right), \left(\frac{7}{9},\frac{1}{9}\right),$ & $9a_2$ \\
& $\left(\frac{1}{3},\frac{1}{9}\right), \left(\frac{2}{9},\frac{1}{3}\right), \left(\frac{8}{9},\frac{2}{9}\right), \left(\frac{2}{3},\frac{8}{9}\right), \left(\frac{7}{9},\frac{2}{3}\right), \left(\frac{1}{9},\frac{7}{9}\right)$ & \\
\hline $(p+q,1-q)$ & $\left(\frac{2}{9},\frac{5}{9}\right), \left(\frac{5}{9},\frac{1}{3}\right), \left(\frac{1}{3},\frac{7}{9}\right), \left(\frac{7}{9},\frac{4}{9}\right), \left(\frac{4}{9},\frac{2}{3}\right), \left(\frac{2}{3},\frac{2}{9}\right),$ & $9a_3$ \\
& $\left(\frac{5}{9},\frac{2}{9}\right), \left(\frac{1}{3},\frac{5}{9}\right), \left(\frac{7}{9},\frac{1}{3}\right), \left(\frac{5}{9},\frac{7}{9}\right), \left(\frac{2}{3},\frac{4}{9}\right), \left(\frac{2}{9},\frac{2}{3}\right)$ & \\
\hline $(0,-1)$ & $\left(\frac{1}{12},\frac{5}{12}\right), \left(\frac{5}{12},\frac{1}{3}\right), \left(\frac{1}{3},\frac{11}{12}\right), \left(\frac{11}{12},\frac{7}{12}\right), \left(\frac{7}{12},\frac{2}{3}\right), \left(\frac{2}{3},\frac{1}{12}\right),$ & 48 \\
& $\left(\frac{5}{12},\frac{1}{12}\right), \left(\frac{1}{3},\frac{5}{12}\right), \left(\frac{11}{12},\frac{1}{3}\right), \left(\frac{7}{12},\frac{11}{12}\right), \left(\frac{2}{3},\frac{7}{12}\right), \left(\frac{1}{12},\frac{2}{3}\right)$ & \\
\hline $\left(\frac{1+\sqrt{13}}{2},1\right)$ & $\left(\frac{1}{13},\frac{4}{13}\right), \left(\frac{4}{13},\frac{3}{13}\right), \left(\frac{3}{13},\frac{12}{13}\right), \left(\frac{12}{13},\frac{9}{13}\right), \left(\frac{9}{13},\frac{10}{13}\right), \left(\frac{10}{13},\frac{1}{13}\right),$ & $13$ \\
& $\left(\frac{4}{13},\frac{1}{13}\right), \left(\frac{3}{13},\frac{4}{13}\right), \left(\frac{12}{13},\frac{3}{13}\right), \left(\frac{9}{13},\frac{12}{13}\right), \left(\frac{10}{13},\frac{9}{13}\right), \left(\frac{1}{13},\frac{10}{13}\right)$ & \\
\hline $\left(\frac{1-\sqrt{13}}{2},1\right)$ & $\left(\frac{2}{13},\frac{7}{13}\right), \left(\frac{7}{13},\frac{5}{13}\right), \left(\frac{5}{13},\frac{11}{13}\right), \left(\frac{11}{13},\frac{6}{13}\right), \left(\frac{6}{13},\frac{8}{13}\right), \left(\frac{8}{13},\frac{2}{13}\right),$ & $13$ \\
& $\left(\frac{7}{13},\frac{2}{13}\right), \left(\frac{5}{13},\frac{7}{13}\right), \left(\frac{11}{13},\frac{5}{13}\right), \left(\frac{6}{13},\frac{11}{13}\right), \left(\frac{8}{13},\frac{6}{13}\right), \left(\frac{2}{13},\frac{8}{13}\right)$ & \\
\hline
\end{tabular}
\caption{$(x,y) \in \mathfrak{D}$ and (the orbits of) the corresponding points $(\theta_1,\theta_2) \in [0,1]^2$. Here $p=2\cos(4\pi/9)$, $q=2\cos(8\pi/9)$, $a_1 = 3 + 2\cos(\pi/9) + 2\cos(2\pi/9)$, $a_2 = 3 - 2\cos(2\pi/9) + 2\sin(\pi/18)$, $a_3 = 3 - 2\cos(\pi/9) - 2\sin(\pi/18)$} \label{Table:subgroupsG2-orbits(theta1,theta2)}
\end{center}
\end{table}

\renewcommand{\arraystretch}{1}

Let $\Omega(\theta_1,\theta_2) := \Phi_1(e^{2\pi i \theta_1},e^{2\pi i \theta_2})^m \Phi_2(e^{2\pi i \theta_1},e^{2\pi i \theta_2})^n \in \mathfrak{D}$ and $\Omega^W(\theta_1,\theta_2)$ its orbit under $W=D_{12}$ given by $\Omega^W(\theta_1,\theta_2) := \sum_{g \in D_{12}} \Omega(g(\theta_1,\theta_2))/12$. We now determine the spectral measures associated with the pairs $(\theta_1,\theta_2)$ that appear in Table \ref{Table:subgroupsG2-orbits(theta1,theta2)}.

We first note that
\begin{align}
\Omega^W(0,0) &= \int_{\bbT^2} \Omega(\theta_1,\theta_2) \, \mathrm{d}_1 \times \mathrm{d}_1, \\
3\Omega^W(0,1/2) &= \int_{\bbT^2} \Omega(\theta_1,\theta_2) \, (4 \, \mathrm{d}_2 \times \mathrm{d}_2 - \mathrm{d}_1 \times \mathrm{d}_1), \\
6\Omega^W(0,1/3) &= \int_{\bbT^2} \Omega(\theta_1,\theta_2) \, (9 \, \mathrm{d}_3 \times \mathrm{d}_3 - 3 \, \mathrm{d}^{(1)}), \\
2\Omega^W(1/3,2/3) &= \int_{\bbT^2} \Omega(\theta_1,\theta_2) \, (3 \mathrm{d}^{(1)} - \mathrm{d}_1 \times \mathrm{d}_1).
\end{align}

\begin{figure}[tb]
\begin{center}
  \includegraphics[width=55mm]{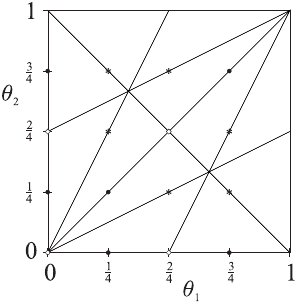}\\
 \caption{The orbits of $(0,1/4) \, \bullet$ and $(1/4,1/2) \, \ast$.} \label{Fig-orbit4}
\end{center}
\end{figure}

The orbits of $(0,1/4)$, $(1/4,1/2)$ are illustrated in Figure \ref{Fig-orbit4}, represented by $\bullet$, $\ast$ respectively. Both orbits lie on the boundary of the fundamental domains of $\bbT^2/D_{12}$, however, only the orbit of $(1/4,1/2)$ lies on the boundary of the fundamental domains of $\bbT^2/S_3$, illustrated in Figure \ref{fig:fund_domain-A2inT2}. Then $6K(0,1/4)\,\Omega^W(0,1/4) = 16 \int_{\bbT^2} \Omega(\theta_1,\theta_2) K(\theta_1,\theta_2) \, \mathrm{d}_4 \times \mathrm{d}_4$ where $K$ is an $S_3$-invariant function on $\bbT^2$ which is zero on the boundaries of the fundamental domains of $\bbT^2/S_3$. Such an $S_3$-invariant function is given by the square $\widetilde{J}^2$ of the Jacobian which appeared for the $A_2$ spectral measures in \cite{evans/pugh:2009v, evans/pugh:2010i}, which is given by $\widetilde{J}(\theta_1,\theta_2) = 4\pi^2(\sin(2\pi(\theta_1+\theta_2))-\sin(2\pi(2\theta_1-\theta_2))-\sin(2\pi(2\theta_2-\theta_1)))$. Now $|\widetilde{J}(0,1/4)| = 8\pi^2$, thus we take
\begin{equation}
K(\theta_1,\theta_2) = \frac{16\widetilde{J}^2(\theta_1,\theta_2)}{64\pi^4} = \frac{\widetilde{J}^2(\theta_1,\theta_2)}{4\pi^4},
\end{equation}
and we have $K(0,1/4) = 16$, so that
\begin{equation}
6\Omega^W(0,1/4) = \int_{\bbT^2} \Omega(\theta_1,\theta_2) K(\theta_1,\theta_2) \, \mathrm{d}_4 \times \mathrm{d}_4.
\end{equation}
We can also obtain an expression for $\Omega^W(1/4,1/2)$, where from Figure \ref{Fig-orbit4} we see that
\begin{equation}
6\Omega^W(1/4,1/2) = \int_{\bbT^2} \Omega(\theta_1,\theta_2) \left( (16 - K(\theta_1,\theta_2)) \, \mathrm{d}_4 \times \mathrm{d}_4 - 4 \, \mathrm{d}_2 \times \mathrm{d}_2 \right).
\end{equation}

\begin{figure}[tb]
\begin{center}
  \includegraphics[width=55mm]{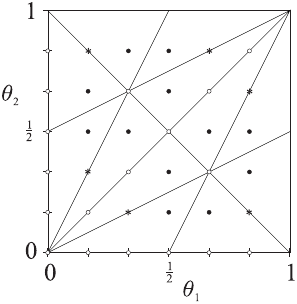}\\
 \caption{The orbits of $(1/6,1/2) \, \bullet $ $(1/6,1/3) \, \ast$.} \label{Fig-orbit6}
\end{center}
\end{figure}

Now
\begin{equation}
12\Omega^W(1/6,1/2) = \int_{\bbT^2} \Omega(\theta_1,\theta_2) \frac{J^2(\theta_1,\theta_2)}{16\pi^4} \, \mathrm{d}_6 \times \mathrm{d}_6,
\end{equation}
as illustrated in Figure \ref{Fig-orbit6} since the Jacobian $J=0$ along the boundaries of the orbit of the fundamental domain, whilst $J^2(g(1/6,1/2)/64\pi^4) = 9$ for all $g \in D_{12}$.
The points $\circ$ in Figure \ref{Fig-orbit6} give the measure $3 \, \mathrm{d}^{(1)}$ whilst the points $\diamond$ in Figure \ref{Fig-orbit4} give the measure $4 \, \mathrm{d}_2 \times \mathrm{d}_2 - \mathrm{d}_1 \times \mathrm{d}_1$, thus we see that
\begin{equation}
6\Omega^W(1/6,1/3) = \int_{\bbT^2} \Omega(\theta_1,\theta_2) \left( 12 \, \mathrm{d}^{(2)} - 3 \, \mathrm{d}^{(1)} - 4 \, \mathrm{d}_2 \times \mathrm{d}_2 + \mathrm{d}_1 \times \mathrm{d}_1 \right).
\end{equation}

\begin{figure}[tb]
\begin{center}
  \includegraphics[width=55mm]{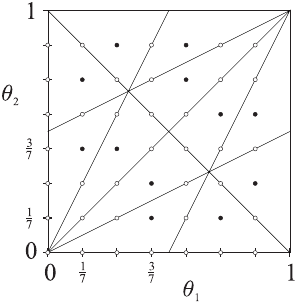}\\
 \caption{The orbit of $(1/7,3/7)$.} \label{Fig-orbit7}
\end{center}
\end{figure}

We have
\begin{equation}
12\Omega^W(1/7,3/7) = \int_{\bbT^2} \Omega(\theta_1,\theta_2) \frac{J^2(\theta_1,\theta_2)}{16\pi^4} \, \mathrm{d}_7 \times \mathrm{d}_7,
\end{equation}
as illustrated in Figure \ref{Fig-orbit7} since the Jacobian $J=0$ along the boundaries of the orbit of the fundamental domain, whilst $J^2(g(1/7,3/7)/64\pi^4) = 49/4$ for all $g \in D_{12}$.

\begin{figure}[tb]
\begin{center}
  \includegraphics[width=55mm]{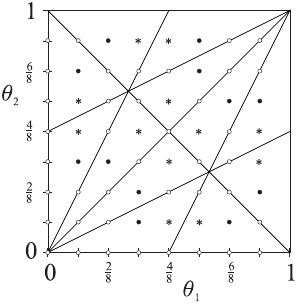}\\
 \caption{The orbits of $(1/8,3/8) \, \bullet$ and $(1/8,1/2) \, \ast$.} \label{Fig-orbit8}
\end{center}
\end{figure}

The orbits of $(1/8,3/8)$, $(1/8,1/2)$ are illustrated in Figure \ref{Fig-orbit8}, represented by $\bullet$, $\ast$ respectively.
We have
\begin{align}
12\Omega^W(1/8,3/8) &= \int_{\bbT^2} \Omega(\theta_1,\theta_2) \left( \sum_{g \in D_{12}} \delta_{g(e^{\pi i/4},e^{3\pi i/4})} \right), \\
12\Omega^W(1/8,1/2) &= \int_{\bbT^2} \Omega(\theta_1,\theta_2) \left( \sum_{g \in D_{12}} \delta_{g(e^{\pi i/4},-1)} \right),
\end{align}
where neither orbit gives a linear combination of the measures in Definition \ref{def:4measures}.

\begin{figure}[tb]
\begin{center}
  \includegraphics[width=55mm]{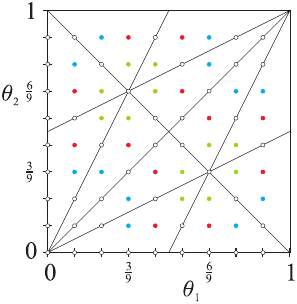}\\
 \caption{The orbits of $(1/9,4/9)\textcolor{red}{\scriptscriptstyle{\bullet}}$, $(1/9,1/3)\textcolor{blue}{\scriptscriptstyle{\bullet}}$, $(2/9,5/9)\textcolor{green}{\scriptscriptstyle{\bullet}}$.} \label{Fig-orbit9}
\end{center}
\end{figure}

Now
\begin{align}
12a_1 \Omega^W(1/9,4/9) &= 2 \int_{\bbT^2} \Omega(\theta_1,\theta_2) \frac{J^2(\theta_1,\theta_2)}{16\pi^4} \, \mathrm{d}^{((9/2))} \\
12a_2 \Omega^W(1/9,1/3) &= 2 \int_{\bbT^2} \Omega(\theta_1,\theta_2) \frac{J^2(\theta_1,\theta_2)}{16\pi^4} \, \mathrm{d}^{((9/4))} \\
12a_3 \Omega^W(2/9,5/9) &= 2 \int_{\bbT^2} \Omega(\theta_1,\theta_2) \frac{J^2(\theta_1,\theta_2)}{16\pi^4} \, \mathrm{d}^{((9))},
\end{align}
as illustrated in Figure \ref{Fig-orbit9} since the Jacobian $J=0$ along the boundaries of the orbit of the fundamental domain.

\begin{figure}[tb]
\begin{center}
  \includegraphics[width=55mm]{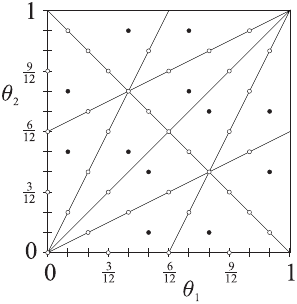}\\
 \caption{The orbit of $(1/12,5/12)$.} \label{Fig-orbit12}
\end{center}
\end{figure}

We have
\begin{equation}
12\Omega^W(1/12,5/12) = \int_{\bbT^2} \Omega(\theta_1,\theta_2) \frac{J^2(\theta_1,\theta_2)}{16\pi^4} \, \mathrm{d}^{(4)},
\end{equation}
as illustrated in Figure \ref{Fig-orbit12} since the Jacobian $J=0$ along the boundaries of the orbit of the fundamental domain, whilst $J^2(g(1/12,5/12)/16\pi^4) = 48$ for all $g \in D_{12}$.

\begin{figure}[tb]
\begin{center}
  \includegraphics[width=55mm]{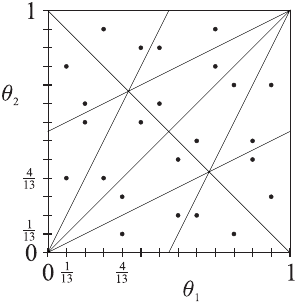}\\
 \caption{The orbits of $(1/13,4/13)$ and $(2/13,7/13)$.} \label{Fig-orbit13}
\end{center}
\end{figure}

Finally, the orbit of the points $(1/13,4/13)$, $(2/13,7/13)$, illustrated in Figure \ref{Fig-orbit13}, do not give a linear combination of the measures in Definition \ref{def:4measures}, thus we have
\begin{align}
12\Omega^W(1/13,4/13) &= \int_{\bbT^2} \Omega(\theta_1,\theta_2) \left( \sum_{g \in D_{12}} \delta_{g(\zeta,\zeta^4)} \right), \\
12\Omega^W(2/13,7/13) &= \int_{\bbT^2} \Omega(\theta_1,\theta_2) \left( \sum_{g \in D_{12}} \delta_{g(\zeta^2,\zeta^7)} \right).
\end{align}

\section{Explicit results} \label{sect:results}
\subsection{Group $PSL(2;7) \cong GL(3;2) \cong \Sigma (168)$} \label{sect:II1}

The subgroup $PSL(2;7)$ of $G_2$ is an irreducible imprimitive group of order 168 which is isomorphic to the group $GL(3;2)$, and also to the subgroup $\Sigma(168)$ of $SU(3)$ which was considered in \cite{evans/pugh:2010i}.

The group $PSL(2;7)$ has irreducible real representations $\Sigma_d$ of dimensions $d = 1,6,7,8$, and two complex conjugate irreducible representations $\Sigma_3, \Sigma_3^{\ast}$ of dimension 3. Its character table is given in Table \ref{table:Character_table-II1} \cite{littlewood:1934}.

\begin{table}[tb]
\begin{center}
\begin{tabular}{|c||c|c|c|c|c|c|} \hline
$C$ & $C_1$ & $C_2$ & $(C_7,C_7^2,C_7^4)$ & $(C_7^3,C_7^5,C_7^6)$ & $(C_4,C_4^3)$ & $(C_3,C_3^2)$ \\
\hline $|C|$ & 1 & 21 & 24 & 24 & 42 & 56 \\
\hline \hline $\Sigma_1$ & 1 & 1 & 1 & 1 & 1 & 1 \\
\hline $\Sigma_3$ & 3 & -1 & $w$ & $\overline{w}$ & 1 & 0 \\
\hline $\Sigma_3^{\ast}$ & 3 & -1 & $\overline{w}$ & $w$ & 1 & 0 \\
\hline $\Sigma_6$ & 6 & 2 & -1 & -1 & 0 & 0  \\
\hline $\Sigma_7$ & 7 & -1 & 0 & 0 & -1 & 1  \\
\hline $\Sigma_8$ & 8 & 0 & 1 & 1 & 0 & -1  \\
\hline
\end{tabular} \\
\caption{Character table for group $PSL(2;7)$, where $w = \eta + \eta^2 + \eta^4 = (-1+i\sqrt{7})/2$, $\eta = e^{2\pi i/7}$.} \label{table:Character_table-II1}
\end{center}
\end{table}

\begin{figure}[tb]
\begin{center}
  \includegraphics[width=110mm]{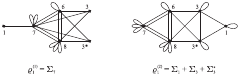}\\
 \caption{The McKay graphs $\mathcal{G}^{\varrho_1^{(i)}}_{PSL(2;7)}$, $i=1,2$.} \label{Fig-McKay_Graph-II1-rho1}
\end{center}
\end{figure}

There are two non-conjugate embeddings of $PSL(2;7)$ in $G_2$ \cite{king/toumazet/wybourne:1999}, given by $\varrho_1^{(1)} = \Sigma_7$ and $\varrho_1^{(2)} = \Sigma_1 + \Sigma_3 + \Sigma_3^{\ast}$.
The McKay graph $\mathcal{G}^{\varrho_1^{(1)}}_{PSL(2;7)}$ for $\varrho_1^{(1)}$ is given in \cite[Figure 1]{he:2003}. We reproduce it in Figure \ref{Fig-McKay_Graph-II1-rho1} for completeness, along with the McKay graph $\mathcal{G}^{\varrho_1^{(2)}}_{PSL(2;7)}$ for $\varrho_1^{(2)}$. We use the notation $n$, $n^{\ast}$ to label the vertices corresponding to the irreducible representations $\Sigma_n$, $\Sigma_n^{\ast}$ respectively.

The eigenvalues of the representation matrices are given in \cite[Tables 4a,b]{king/toumazet/wybourne:1999}.
The decomposition of the Kronecker square of $\varrho_1^{(i)}$ into irreducibles is given by
$$(\varrho_1^{(1)})^2 = \mathrm{id} + \varrho_1^{(1)} + \Sigma_3 + \Sigma_3^{\ast} + 2\Sigma_6 + \Sigma_7 + 2 \Sigma_8, \qquad
(\varrho_1^{(2)})^2 = \mathrm{id} + \varrho_1^{(2)} + \Sigma_1 + 2\Sigma_3 + 2\Sigma_3^{\ast} + 2\Sigma_6 + 2 \Sigma_8,$$
where $\mathrm{id} = \Sigma_1$.
From dimension considerations there are thus two candidates for the fourteen-dimensional representation $\varrho_2^{(i)}$, which are given by $\Sigma_3 + \Sigma_3^{\ast} + \Sigma_8$ and $\Sigma_6 + \Sigma_8$ for both $i=1,2$.
However, as discussed in Section \ref{sect:subgroupsG2}, since $\chi_{\varrho_2}(C) = \Phi_2(t_1^C,t_2^C)$, where $(t_1^{C},t_2^{C})$ is a pair from the set of eigenvalues of group elements from the conjugacy class $C$, from knowledge of the eigenvalues from \cite{king/toumazet/wybourne:1999} we see that the decomposition of the fundamental fourteen-dimensional representation into irreducible representations of $PSL(2;7)$ is given by $\varrho_2 := \varrho_2^{(i)} = \Sigma_3 + \Sigma_3^{\ast} + \Sigma_8$ for both $i=1,2$.
The values of $x^{(i)} = \chi_{\varrho_1^{(i)}}(C) \in [-2,7]$, $y = \chi_{\rho_2}(C) \in [-2,14]$ for $PSL(2;7)$ are given in Table \ref{table:(x,y)-II1}, along with the values of $J^2/64\pi^4$ for the corresponding pairs $(\theta_1,\theta_2) \in [0,1]^2$ obtained from Table \ref{Table:subgroupsG2-orbits(theta1,theta2)}.

\begin{table}[tb]
\begin{center}
\begin{tabular}{|c||c|c|c|c|c|c|} \hline
$C$ & $C_1$ & $C_2$ & $(C_7,C_7^2,C_7^4)$ & $(C_7^3,C_7^5,C_7^6)$ & $(C_4,C_4^3)$ & $(C_3,C_3^2)$ \\
\hline $\chi_{\varrho_1^{(1)}}(C) \in [-2,7]$ & 7 & -1 & 0 & 0 & -1 & 1  \\
\hline $\chi_{\varrho_1^{(2)}}(C) \in [-2,7]$ & 7 & -1 & 0 & 0 & 3 & 1  \\
\hline $\chi_{\varrho_2}(C) \in [-2,14]$ & 14 & -2 & 0 & 0 & 2 & -1 \\
\hline $J^2(\theta_1,\theta_2)/64\pi^4$ & 0 & 0 & 49/4 & 49/4 & 0 & 0 \\
\hline
\end{tabular} \\
\caption{$\chi_{\varrho_j}(C)$ for group $PSL(2;7)$, $j=1,2$.} \label{table:(x,y)-II1}
\end{center}
\end{table}

From (\ref{eqn:moments-subgroupG2}) and Tables \ref{table:Character_table-II1}, \ref{table:(x,y)-II1} and \ref{Table:subgroupsG2-orbits(theta1,theta2)}, we see that
$$\varsigma_{m,n} = \frac{1}{168} \Omega^W(0,0) + \frac{21}{168} \Omega^W(0,1/2) + \frac{56}{168} \Omega^W(0,1/3) + \frac{42}{168} \Omega' + \frac{24+24}{168} \Omega^W(1/7,3/7),$$
where $\Omega'$ is $\Omega^W(1/4,1/2)$ for $\varrho_1^{(1)}$, and $\Omega^W(0,1/4)$ for $\varrho_1^{(2)}$.
Then we obtain

\begin{Thm} \label{thm:measureII1}
The joint spectral measure (over $\bbT^2$) for the non-conjugate embeddings of the projective special linear group $PSL(2;7)$ over the finite field $\mathbb{F}_7$ into the fundamental representations of $G_2$ is
\begin{equation} \label{eqn:measureII1}
\mathrm{d}\varepsilon = \frac{1}{672\pi^4} J^2 \, \mathrm{d}_7 \times \mathrm{d}_7 + \frac{1}{24} K' \, \mathrm{d}_4 \times \mathrm{d}_4 + \frac{1}{2} \, \mathrm{d}_3 \times \mathrm{d}_3 + \frac{1}{6} \beta \, \mathrm{d}_2 \times \mathrm{d}_2 - \frac{1}{28} \, \mathrm{d}_1 \times \mathrm{d}_1 - \frac{1}{6} \, \mathrm{d}^{(1)}
\end{equation}
where for the embedding of $PSL(2;7)$ in $G_2$ given by $\varrho_1^{(1)} = \Sigma_7$ we have $K' = 16-K$ and $\beta = 1$, whilst for the embedding given by $\varrho_1^{(2)} = \Sigma_1 + \Sigma_3 + \Sigma_3^{\ast}$ we have $K' = K$ and $\beta = 0$, where $K(\theta_1,\theta_2) = (\sin(2\pi(\theta_1+\theta_2))-\sin(2\pi(2\theta_1-\theta_2))-\sin(2\pi(2\theta_2-\theta_1)))^2$, $\mathrm{d}_m$ is the uniform measure over $m^{\mathrm{th}}$ roots of unity and $\mathrm{d}^{(k+4)}$ is the uniform measure on the points in $C_k^W$.
\end{Thm}

\begin{Rem}
Note that measure in Theorem \ref{thm:measureII1} for the second embedding $\varrho_1^{(2)}$ of $PSL(2;7)$ in $G_2$ is precisely that for $\Sigma (168) \subset SU(3)$ given in \cite[Theorem 16]{evans/pugh:2010i}. However, (\ref{eqn:measureII1}) has a cleaner expression than that given in \cite{evans/pugh:2010i}, because here we were able to use the Jacobian $J$ for $G_2$ which is also 0 along the diagonal, whereas the Jacobian for $SU(3)$ (essentially $K$ in Theorem \ref{thm:measureII1}) is non-zero along the diagonal. \end{Rem}

\subsection{Group $\bbZ_2^3 \cdot PSL(2;7)$} \label{sect:II2}

The subgroup $\bbZ_2^3 \cdot PSL(2;7)$ of $G_2$ is a nonsplit extension of $\bbZ_2^3$ by $PSL(2;7)$. It is an irreducible imprimitive group of order 1344.
It has eleven irreducible representations (nine real and two complex conjugate representations). As a subgroup of $GL(7)$ it is generated by the permutations $(1234567)$ and $(124)(365)$, and the signed permutation  $\delta^-_{\{1,2,4,7\}}(12)(36)$, where for a set $S$ the notation $\delta^-_S$ means that the entries in the rows indexed by elements of $S$ are multiplied by $-1$ \cite{coxeter:1946, cohen/wales:1983}.
Its character table, generated by the GAP computational discrete algebra system, is given in Table \ref{table:Character_table-II2} (see also \cite{he:2003})
where elements in $C_4$, $C_4^{\prime}$, $C_8$, $C_8^{\prime}$, $C_7$, $C_7^{\prime}$, $C_6$, $C_3$ are of cycle type $(C_4,C_4^3)$, $(C_4^{\prime},C_4^{\prime3})$, $(C_8,C_8^3,C_8^5,C_8^7)$, $(C_8^{\prime},C_8^{\prime3},C_8^{\prime5},C_8^{\prime7})$, $(C_7,C_7^2,C_7^4)$, $(C_7^3,C_7^5,C_7^6)$, $(C_6,C_6^5)$, $(C_3,C_3^2)$ respectively.

\begin{table}[tb]
\begin{center}
\begin{tabular}{|c||c|c|c|c|c|c|c|c|c|c|c|} \hline
$C$                             & $C_1$ & $C_2$ & $C_4$ & $C_4^{\prime}$ & $C_8$ & $C_8^{\prime}$ & $C_2^{\prime}$   & $C_7$            & $C_7^{\prime}$ & $C_6$ & $C_3$    \\
\hline $|C|$                    & 1     & 7     & 42    & 42             & 168   & 168            & 84               & 192              & 192            & 224   & 224      \\
\hline \hline $\Sigma_1$        & 1     & 1     & 1     & 1              & 1     & 1              & 1                & 1                & 1              & 1     & 1        \\
\hline $\Sigma_3$               & 3     & 3     & -1    & -1             & 1     & 1              & -1               & $w$              & $\overline{w}$ & 0     & 0        \\
\hline $\Sigma_3^{\ast}$        & 3     & 3     & -1    & -1             & 1     & 1              & -1               & $\overline{w}$   & $w$            & 0     & 0        \\
\hline $\Sigma_6$               & 6     & 6     & 2     & 2              & 0     & 0              & 2                & -1               & -1             & 0     & 0        \\
\hline $\Sigma_7^{(1)}$         & 7     & -1    & -1    & 3              & 1     & -1             & -1               & 0                & 0              & -1    & 1        \\
\hline $\Sigma_7^{(1)\prime}$   & 7     & -1    & 3     & -1             & -1    & 1              & -1               & 0                & 0              & -1    & 1        \\
\hline $\Sigma_7^{(2)}$         & 7     & 7     & -1    & -1             & -1    & -1             & -1               & 0                & 0              & 1     & 1        \\
\hline $\Sigma_8$               & 8     & 8     & 0     & 0              & 0     & 0              & 0                & 1                & 1              & -1    & -1       \\
\hline $\Sigma_{14}$            & 14    & -2    & 2     & 2              & 0     & 0              & -2               & 0                & 0              & 1     & -1       \\
\hline $\Sigma_{21}$            & 21    & -3    & -3    & 1              & -1    & 1              & 1                & 0                & 0              & 0     & 0        \\
\hline $\Sigma_{21}'$           & 21    & -3    & 1     & -3             & 1     & -1             & 1                & 0                & 0              & 0     & 0        \\
\hline
\end{tabular} \\
\caption{Character table for $\bbZ_2^3 \cdot PSL(2;7)$, where $w = \eta + \eta^2 + \eta^4 = (-1+i\sqrt{7})/2$, $\eta = e^{2\pi i/7}$.} \label{table:Character_table-II2}
\end{center}
\end{table}

There are five non-conjugate seven-dimensional representations, $\Sigma_1 + \Sigma_3 + \Sigma_3^{\ast}$, $\Sigma_1 + \Sigma_6$, $\Sigma_7^{(1)}$, $\Sigma_7^{(1)\prime}$ and $\Sigma_7^{(2)}$, however of these, only $\varrho_1^{(1)} := \Sigma_7^{(1)}$ and $\varrho_1^{(2)} := \Sigma_7^{(1)\prime}$ are faithful.
Both these satisfy the condition that $\varrho_1^{(i)}$ appears in the decomposition of $(\varrho_1^{(i)})^2$, and their eigenvalues are given in Table \ref{table:evalues-II2}.
As described in Section \ref{sect:subgroupsG2}, these eigenvalues can be determined from the character table of $\bbZ_2^3 \cdot PSL(2;7)$. The additional information that is needed is to note that the eigenvalues for group elements in $C_4$ and $C_4^{\prime}$ square to those for elements in $C_2$, those for $C_8$ square to those for $C_4$, those for $C_8^{\prime}$ square to those for $C_4^{\prime}$, whilst those for $C_6$ square to those for $C_3$ and also cube to those for $C_2$.

\begin{table}[tb]
\begin{center}
\begin{tabular}{|c||c|c|} \hline
                                & $\varrho_1^{(1)} = \Sigma_7^{(1)}$            & $\varrho_1^{(2)} = \Sigma_7^{(1)\prime}$      \\ 
\hline \hline   $C_1$           & $(1,1,1,1,1,1,1)$                             & $(1,1,1,1,1,1,1)$                             \\ 
\hline          $C_2$           & $(1,1,1,-1,-1,-1,-1)$                         & $(1,1,1,-1,-1,-1,-1)$                         \\ 
\hline          $C_4$           & $(1,-1,-1,i,i,-i,-i)$                         & $(1,1,1,i,i,-i,-i)$                           \\ 
\hline          $C_4^{\prime}$  & $(1,1,1,i,i,-i,-i)$                           & $(1,-1,-1,i,i,-i,-i)$                         \\ 
\hline          $C_8$           & $(1,-1,-1,\nu,\nu^3,\nu^5,\nu^7)$             & $(1,i,-i,\nu,\nu^3,\nu^5,\nu^7)$              \\ 
\hline          $C_8^{\prime}$  & $(1,i,-i,\nu,\nu^3,\nu^5,\nu^7)$              & $(1,-1,-1,\nu,\nu^3,\nu^5,\nu^7)$             \\ 
\hline          $C_2^{\prime}$  & $(1,1,1,-1,-1,-1,-1)$                         & $(1,1,1,-1,-1,-1,-1)$                         \\ 
\hline          $C_7$           & $(1,\eta,\eta^2,\eta^3,\eta^4,\eta^5,\eta^6)$ & $(1,\eta,\eta^2,\eta^3,\eta^4,\eta^5,\eta^6)$ \\ 
\hline          $C_7^{\prime}$  & $(1,\eta,\eta^2,\eta^3,\eta^4,\eta^5,\eta^6)$ & $(1,\eta,\eta^2,\eta^3,\eta^4,\eta^5,\eta^6)$ \\ 
\hline          $C_6$           & $(1,-1,-1,\omega,-\omega,\omega^2,-\omega^2)$ & $(1,-1,-1,\omega,-\omega,\omega^2,-\omega^2)$ \\ 
\hline          $C_3$           & $(1,1,1,\omega,\omega,\omega^2,\omega^2)$     & $(1,1,1,\omega,\omega,\omega^2,\omega^2)$     \\ 
\hline
\end{tabular} \\
\caption{Eigenvalues of group elements in each conjugacy class of $\bbZ_2^3 \cdot PSL(2;7)$ for the irreducible representations $\varrho_1^{(i)}$, $i=1,2$, where $\omega = e^{2\pi i/3}$, $\nu = e^{2\pi i/8}$ and $\eta = e^{2\pi i/7}$.} \label{table:evalues-II2}
\end{center}
\end{table}

We present one choice of eigenvalues $(t_1^C,t_2^C)$ in Table \ref{table:t1,t2-II2}, along with the characters of the corresponding fourteen-dimensional representation for these eigenvalues. In both cases we see that the character agrees with the irreducible representation $\Sigma_{14}$.

\begin{table}[tb]
\begin{center}
\begin{tabular}{|c||c|c|c|c|c|} \hline
& \multicolumn{2}{|c|}{$\varrho_1^{(1)}$} & \multicolumn{2}{|c|}{$\varrho_1^{(2)}$} & \\
       $C$              & $(t_1^C,t_2^C)$   & $\chi_{\varrho_2^{(1)}}(C)$   & $(t_1^C,t_2^C)$   & $\chi_{\varrho_2^{(2)}}(C)$   & $\chi_{\Sigma_{14}}(C)$ \\
\hline \hline $C_1$     & $(1,1)$           & 14                            & $(1,1)$           & 14                            & 14                         \\
\hline $C_2$            & $(1,-1)$          & -2                            & $(1,-1)$          & -2                            & -2                         \\
\hline $C_4$            & $(-1,i)$          & 2                             & $(1,i)$           & 2                             & 2                          \\
\hline $C_4^{\prime}$   & $(1,i)$           & 2                             & $(-1,i)$          & 2                             & 2                          \\
\hline $C_8$            & $(-1,\nu)$        & 0                             & $(i,\nu)$         & 0                             & 0                          \\
\hline $C_8^{\prime}$   & $(i,\nu)$         & 0                             & $(-1,\nu)$        & 0                             & 0                          \\
\hline $C_2^{\prime}$   & $(1,-1)$          & -2                            & $(1,-1)$          & -2                            & -2                         \\
\hline $C_7$            & $(\eta,\eta^5)$   & 0                             & $(\eta,\eta^5)$   & 0                             & 0                          \\
\hline $C_7^{\prime}$   & $(\eta^2,\eta^3)$ & 0                             & $(\eta^2,\eta^3)$ & 0                             & 0                          \\
\hline $C_6$            & $(-1,\omega)$     & 1                             & $(-1,\omega)$     & 1                             & 1                          \\
\hline $C_3$            & $(1,\omega)$      & -1                            & $(1,\omega)$      & -1                            & -1                         \\
\hline
\end{tabular} \\
\caption{Choice of eigenvalues $(t_1^C,t_2^C)$ for $\varrho_1^{(i)}$, $i=1,2$, and corresponding values of $\chi_{\varrho_2^{(i)}}(C)$.} \label{table:t1,t2-II2}
\end{center}
\end{table}

\begin{figure}[tb]
\begin{center}
  \includegraphics[width=70mm]{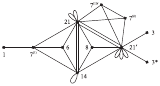}\\
 \caption{The McKay graph $\mathcal{G}^{\varrho_1^{(1)}}_{\bbZ_2^3 \cdot PSL(2;7)}$.} \label{Fig-McKay_Graph-II2-rho1}
\end{center}
\end{figure}

Thus there are precisely two non-conjugate faithful irreducible representations $\varrho_1^{(i)}$, $i=1,2$, giving embeddings of $\bbZ_2^3 \cdot PSL(2;7)$ in $G_2$.
The McKay graph $\mathcal{G}^{\varrho_1^{(1)}}_{PSL(2;7)}$ for $\varrho_1^{(1)}$ is given in \cite[Figure 1]{he:2003}. We reproduce it in Figure \ref{Fig-McKay_Graph-II2-rho1} for completeness. We use the notation $n$, $n^{\ast}$, $n^{(i)\prime}$ to label the vertices corresponding to the irreducible representations $\Sigma_n$, $\Sigma_n^{\ast}$, $\Sigma_n^{(i)\prime}$ respectively. The McKay graph for $\varrho_1^{(2)}$ is similar.

\begin{table}[tb]
\begin{center}
\begin{tabular}{|c||c|c|c|c|c|c|c|c|c|c|c|} \hline
$C$ & $C_1$ & $C_2$ & $C_4$ & $C_4^{\prime}$ & $C_8$ & $C_8^{\prime}$ & $C_2^{\prime}$ & $C_7$ & $C_7^{\prime}$ & $C_6$ & $C_3$ \\
\hline $\chi_{\varrho_1^{(1)}}(\Gamma_j) \in [-2,7]$ & 7  & -1 & -1 & 3  & 1  & -1 & -1 & 0     & 0     & -1 & 1   \\
\hline $\chi_{\varrho_1^{(2)}}(\Gamma_j) \in [-2,7]$ & 7  & -1 & 3  & -1 & -1 & 1  & -1 & 0     & 0     & -1 & 1   \\
\hline $\chi_{\varrho_2}(\Gamma_j) \in [-2,14]$      & 14 & -2 & 2  & 2  & 0  & 0  & -2 & 0     & 0     & 1  & -1  \\
\hline $J^2(\theta_1,\theta_2)/64\pi^4$              & 0  & 0  & 0  & 0  & 8  & 8  & 0  & 49/4  & 49/4  & 9  & 0   \\
\hline
\end{tabular} \\
\caption{$\chi_{\varrho_j}(C)$ for group $\bbZ_2^3 \cdot PSL(2;7)$, $j=1,2$.} \label{table:(x,y)-II2}
\end{center}
\end{table}

The values of $x^{(i)} = \chi_{\varrho_1^{(i)}}(C) \in [-2,7]$, $y = \chi_{\rho_2}(C) \in [-2,14]$ for $\bbZ_2^3 \cdot PSL(2;7)$ are given in Table \ref{table:(x,y)-II1}, along with the values of $J^2/64\pi^4$ for the corresponding pairs $(\theta_1,\theta_2) \in [0,1]^2$ obtained from Table \ref{Table:subgroupsG2-orbits(theta1,theta2)}.
Then from (\ref{eqn:moments-subgroupG2}) and Tables \ref{table:Character_table-II2}, \ref{table:(x,y)-II2} and \ref{Table:subgroupsG2-orbits(theta1,theta2)}, we see that
\begin{align*}
\varsigma_{m,n} & = \frac{1}{1344} \Omega^W(0,0) + \frac{7+84}{1344} \Omega^W(0,1/2) + \frac{42}{1344} \Omega^W(0,1/4) + \frac{42}{1344} \Omega^W(1/4,1/2)  \\
& \quad + \frac{168}{1344} \Omega^W(1/8,1/2) + \frac{168}{1344} \Omega^W(1/8,3/8) + \frac{192+192}{1344} \Omega^W(1/7,3/7) \\
& \quad + \frac{224}{1344} \Omega^W(1/6,1/2) + \frac{224}{1344} \Omega^W(0,1/3) ,
\end{align*}
where $\Omega^W(\theta_1,\theta_2)$ is as in Section \ref{sect:II1}.
Note that $12(\Omega^W(1/8,1/2) + \Omega^W(1/8,3/8)) = 8 \int_{\bbT^2} \Omega(\theta_1,\theta_2) (J^2/64\pi^4) \, \mathrm{d}_8 \times \mathrm{d}_8$, since the Jacobian $J=0$ along the boundaries of the orbit of the fundamental domain, whilst $J^2(g(1/8,1/2))/64\pi^4 = J^2(g(1/8,3/8))/64\pi^4 = 8$ for all $g \in D_{12}$.
Thus we obtain:

\begin{Thm}
The joint spectral measure (over $\bbT^2$) for both non-conjugate embeddings of $\bbZ_2^3 \cdot PSL(2;7)$ into the fundamental representations of $G_2$ are given by
\begin{align} \label{eqn:measureII2}
\mathrm{d}\varepsilon &= \frac{1}{768\pi^4} J^2 \, \mathrm{d}_8 \times \mathrm{d}_8 + \frac{1}{672\pi^4} J^2 \, \mathrm{d}_7 \times \mathrm{d}_7 + \frac{1}{1152\pi^4} J^2 \, \mathrm{d}_6 \times \mathrm{d}_6 + \frac{1}{12} \, \mathrm{d}_4 \times \mathrm{d}_4  \nonumber \\
& \qquad + \frac{1}{4} \beta \, \mathrm{d}_3 \times \mathrm{d}_3 + \frac{5}{72} \, \mathrm{d}_2 \times \mathrm{d}_2 - \frac{11}{504} \, \mathrm{d}_1 \times \mathrm{d}_1 - \frac{1}{12} \, \mathrm{d}^{(1)}
\end{align}
where $\mathrm{d}_m$ is the uniform measure over $m^{\mathrm{th}}$ roots of unity and $\mathrm{d}^{(k+4)}$ is the uniform measure on the points in $C_k^W$..
\end{Thm}

\subsection{Group $PGL(2;7)$}

The subgroup $PGL(2;7)$ of $G_2$ is an irreducible primitive group of order 336.
It has nine irreducible representations, all real, and its character table is given in Table \ref{table:Character_table-IP3} \cite{collins:1990, he:2003}.

\begin{table}[tb]
\begin{center}
\begin{tabular}{|c||c|c|c|c|c|c|c|c|c|} \hline
$C$ & $C_1$ & $C_2$ & $C_2^{\prime}$ & $C_8$ & $C_8^{\prime}$ & $C_4$ & $C_7$ & $C_6$ & $C_3$ \\
\hline $|C|$ & 1 & 21 & 28 & 42 & 42 & 42 & 48 & 56 & 56 \\
\hline \hline $\Sigma_1$ & 1 & 1 & 1 & 1 & 1 & 1 & 1 & 1 & 1 \\
\hline $\Sigma_1^{\prime}$ & 1 & 1 & -1 & -1 & -1 & 1 & 1 & -1 & 1 \\
\hline $\Sigma_6^{(1)}$ & 6 & -2 & 0 & 0 & 0 & 2 & -1 & 0 & 0 \\
\hline $\Sigma_6^{(2)}$ & 6 & 2 & 0 & $\sqrt{2}$ & $-\sqrt{2}$ & 0 & -1 & 0 & 0 \\
\hline $\Sigma_6^{(2)\prime}$ & 6 & 2 & 0 & $-\sqrt{2}$ & $\sqrt{2}$ & 0 & -1 & 0 & 0 \\
\hline $\Sigma_7$ & 7 & -1 & 1 & -1 & -1 & -1 & 0 & 1 & 1 \\
\hline $\Sigma_7^{\prime}$ & 7 & -1 & -1 & 1 & 1 & -1 & 0 & -1 & 1 \\
\hline $\Sigma_8$ & 8 & 0 & 2 & 0 & 0 & 0 & 1 & -1 & -1 \\
\hline $\Sigma_8^{\prime}$ & 8 & 0 & -2 & 0 & 0 & 0 & 1 & 1 & -1 \\
\hline
\end{tabular} \\
\caption{Character table for $PGL(2;7)$.} \label{table:Character_table-IP3}
\end{center}
\end{table}

There are eight non-conjugate seven-dimensional representations, $\gamma_1^{(1)} = \Sigma_1 + \Sigma_6^{(1)}$, $\gamma_1^{(2)} = \Sigma_1 + \Sigma_6^{(2)}$, $\gamma_1^{(3)} = \Sigma_1 + \Sigma_6^{(2)\prime}$, $\gamma_1^{(4)} = \Sigma_1^{\prime} + \Sigma_6^{(1)}$, $\gamma_1^{(5)} = \Sigma_1^{\prime} + \Sigma_6^{(2)}$, $\gamma_1^{(6)} = \Sigma_1^{\prime} + \Sigma_6^{(2)\prime}$, $\gamma_1^{(7)} = \Sigma_7$ and $\gamma_1^{(8)} = \Sigma_7^{\prime}$.
The decomposition of the Kronecker squares of the $\gamma_1^{(i)}$ are given by:
\begin{align*}
(\gamma_1^{(1)})^2 & = \mathrm{id} + \gamma_1^{(1)} + \Sigma_1^{\prime} + 2\Sigma_6^{(1)} + \Sigma_6^{(2)} + \Sigma_6^{(2)\prime} + \Sigma_8 + \Sigma_8^{\prime}, \\
(\gamma_1^{(2)})^2 & = \mathrm{id} + \gamma_1^{(2)} + 2\Sigma_6^{(2)} + \Sigma_6^{(2)\prime} + \Sigma_7^{\prime} + \Sigma_8 + \Sigma_8^{\prime}, \\
(\gamma_1^{(4)})^2 & = \mathrm{id} + \gamma_1^{(4)} + \Sigma_1 + 2\Sigma_6^{(1)} + \Sigma_6^{(2)} + \Sigma_6^{(2)\prime} + \Sigma_8 + \Sigma_8^{\prime}, \\
(\gamma_1^{(5)})^2 & = \mathrm{id} + \Sigma_1 + 3\Sigma_6^{(2)} + \Sigma_6^{(2)\prime} + \Sigma_7^{\prime} + \Sigma_8 + \Sigma_8^{\prime}, \\
(\gamma_1^{(7)})^2 & = \mathrm{id} + \gamma_1^{(7)} + \Sigma_1 + \Sigma_6^{(1)} + \Sigma_6^{(2)} + \Sigma_6^{(2)\prime} + \Sigma_7^{\prime} + \Sigma_8 + \Sigma_8^{\prime}, \\
(\gamma_1^{(8)})^2 & = \mathrm{id} + \gamma_1^{(8)} + \Sigma_1 + \Sigma_6^{(1)} + \Sigma_6^{(2)} + \Sigma_6^{(2)\prime} + \Sigma_7 + \Sigma_8 + \Sigma_8^{\prime},
\end{align*}
where $\mathrm{id} = \Sigma_1$.
Note, we have omitted the decompositions for $\gamma_1^{(3)}$, $\gamma_1^{(6)}$, since from the character table we see that these representations are essentially the same as $\gamma_1^{(2)}$, $\gamma_1^{(5)}$ respectively.
Then we see that $\gamma_1^{(i)}$ does not appear in the decomposition of $(\gamma_1^{(i)})^2$ for $i=5,6$, therefore they are not embeddings of $PGL(2;7)$ in $G_2$.
From dimension considerations, we see that candidates $\gamma_2^{(j)}$ for $\varrho_2$ are $\gamma_2^{(1)}=\Sigma_6^{(1)}+\Sigma_8$, $\gamma_2^{(2)}=\Sigma_6^{(1)}+\Sigma_8^{\prime}$, $\gamma_2^{(3)}=\Sigma_6^{(2)}+\Sigma_8$, $\gamma_2^{(4)}=\Sigma_6^{(2)}+\Sigma_8^{\prime}$, $\gamma_2^{(5)}=\Sigma_6^{(2)\prime}+\Sigma_8$ and $\gamma_2^{(6)}=\Sigma_6^{(2)\prime}+\Sigma_8^{\prime}$, where $1 \leq j \leq 6$ for $\gamma_1^{(i)}$ when $i=1,4,7,8$, and $3 \leq j \leq 6$ when $i=2$.
Then with $(x_i^C,y_j^C) = (\chi_{\gamma_1^{(i)}}(C),\chi_{\gamma_2^{(j)}}(C))$, we see that $(x_i^C,y_j^C) \not \in \mathfrak{D}$ for any candidate $\gamma_2^{(j)}$ in the cases $i=1,2,3,7$ when $C=C_2^{\prime}$. Thus $\gamma_1^{(i)}$ cannot define an embedding of $PGL(2;7)$ in $G_2$ for $i=1,2,3,7$.
We also see that $(x_8^C,y_j^C) \not \in \mathfrak{D}$ for $j=3,4,5,6$.

Thus we have candidates $(\gamma_1^{(i)},\gamma_2^{(j)})$ for $(\varrho_1,\varrho_2)$ when $i=4,8$, where $1 \leq j \leq 6$ for $i=4$ and $j \in \{ 1,2 \}$ for $i=8$.

\begin{table}[tb]
\begin{center}
\begin{tabular}{|c||c|c|c|c|c|} \hline
 & $\Sigma_1$ & $\Sigma_1^{\prime}$ & $\Sigma_6^{(1)}$ & $\Sigma_6^{(2)}$ & $\Sigma_7$ \\
\hline \hline $C_1$ & 1 & 1 & $(1,1,1,1,1,1)$ & $(1,1,1,1,1,1)$ & $(1,1,1,1,1,1,1)$ \\
\hline $C_2$ & 1 & 1 & $(1,1,-1,-1,-1,-1)$ & $(1,1,1,1,-1,-1)$ & $(1,1,1,-1,-1,-1,-1)$ \\
\hline $C_2^{\prime}$ & 1 & -1 & $(1,1,1,-1,-1,-1)$ & $(1,1,1,-1,-1,-1)$ & $(1,1,1,1,-1,-1,-1)$ \\
\hline $C_8$ & 1 & -1 & $(1,-1,\nu,\nu^3,\nu^5,\nu^7)$ & $(1,-1,i,-i,\nu,\nu^7)$ & $(-1,i,-i,\nu,\nu^3,\nu^5,\nu^7)$ \\
\hline $C_8^{\prime}$ & 1 & -1 & $(1,-1,\nu,\nu^3,\nu^5,\nu^7)$ & $(1,-1,i,-i,\nu^3,\nu^5)$ & $(-1,i,-i,\nu,\nu^3,\nu^5,\nu^7)$ \\
\hline $C_4$ & 1 & 1 & $(1,1,i,i,-i,-i)$ & $(1,1,-1,-1,i,-i)$ & $(1,-1,-1,i,i,-i,-i)$ \\
\hline $C_7$ & 1 & 1 & $(\eta,\eta^2,\eta^3,\eta^4,\eta^5,\eta^6)$ & $(\eta,\eta^2,\eta^3,\eta^4,\eta^5,\eta^6)$ & $(1,\eta,\eta^2,\eta^3,\eta^4,\eta^5,\eta^6)$ \\
\hline $C_6$ & 1 & -1 & $(1,-1,\mu,\mu^2,\mu^4,\mu^5)$ & $(1,-1,\mu,\mu^2,\mu^4,\mu^5)$ & $(1,1,-1,\mu,\mu^2,\mu^4,\mu^5)$ \\
\hline $C_3$ & 1 & 1 & $(1,1,\omega,\omega,\omega^2,\omega^2)$ & $(1,1,\omega,\omega,\omega^2,\omega^2)$ & $(1,1,1,\omega,\omega,\omega^2,\omega^2)$ \\
\hline
\end{tabular} \\
\begin{tabular}{|c||c|c|c|} \hline
 & $\Sigma_7^{\prime}$ & $\Sigma_8$ & $\Sigma_8^{\prime}$ \\
\hline \hline $C_1$ & $(1,1,1,1,1,1,1)$ & $(1,1,1,1,1,1,1,1)$ & $(1,1,1,1,1,1,1,1)$ \\
\hline $C_2$ & $(1,1,1,-1,-1,-1,-1)$ & $(1,1,1,1,-1,-1,-1,-1)$ & $(1,1,1,1,-1,-1,-1,-1)$ \\
\hline $C_2^{\prime}$ & $(1,1,1,-1,-1,-1,-1)$ & $(1,1,1,1,1,-1,-1,-1)$ & $(1,1,1,-1,-1,-1,-1,-1)$ \\
\hline $C_8$ & $(1,i,-i,\nu,\nu^3,\nu^5,\nu^7)$ & $(1,-1,i,-i,\nu,\nu^3,\nu^5,\nu^7)$ & $(1,-1,i,-i,\nu,\nu^3,\nu^5,\nu^7)$ \\
\hline $C_8^{\prime}$ & $(1,i,-i,\nu,\nu^3,\nu^5,\nu^7)$ & $(1,-1,i,-i,\nu,\nu^3,\nu^5,\nu^7)$ & $(1,-1,i,-i,\nu,\nu^3,\nu^5,\nu^7)$ \\
\hline $C_4$ & $(1,-1,-1,i,i,-i,-i)$ & $(1,1,-1,-1,i,i,-i,-i)$ & $(1,1,-1,-1,i,i,-i,-i)$ \\
\hline $C_7$ & $(1,\eta,\eta^2,\eta^3,\eta^4,\eta^5,\eta^6)$ & $(1,1,\eta,\eta^2,\eta^3,\eta^4,\eta^5,\eta^6)$ & $(1,1,\eta,\eta^2,\eta^3,\eta^4,\eta^5,\eta^6)$ \\
\hline $C_6$ & $(1,-1,-1,\mu,\mu^2,\mu^4,\mu^5)$ & $(1,-1,\mu,\mu^2,\mu^2,\mu^4,\mu^4,\mu^5)$ & $(1,-1,\mu,\mu,\mu^2,\mu^4,\mu^5,\mu^5)$ \\
\hline $C_3$ & $(1,1,1,\omega,\omega,\omega^2,\omega^2)$ & $(1,1,\omega,\omega,\omega,\omega^2,\omega^2,\omega^2)$ & $(1,1,\omega,\omega,\omega,\omega^2,\omega^2,\omega^2)$ \\
\hline
\end{tabular} \\
\caption{Eigenvalues of group elements in each conjugacy class of $PGL(2;7)$, where $\omega = e^{2\pi i/3}$, $\mu = e^{2\pi i/6}$ and $\eta = e^{2\pi i/7}$.} \label{table:evalues-IP3}
\end{center}
\end{table}

We now consider the eigenvalues of the representation matrices to determine which of the remaining $\gamma_1^{(i)}$ are embeddings of $PSL(2;7)$ in $G_2$.
These eigenvalues are given in Table \ref{table:evalues-IP3}. As described in Section \ref{sect:subgroupsG2}, these eigenvalues can be determined from the character table of $PGL(2;7)$. The additional information that is needed is to note that the eigenvalues for group elements in $C_4$ square to those for elements in $C_2$, those for $C_8$, $C_8^{\prime}$ square to those for $C_4$, those for $C_6$ square to those for $C_3)$ and also cube to those for $C_2^{\prime}$.
These observations follow from the fact that, for example, it is impossible to choose eigenvalues for group elements in $C_4$ which square to those for elements in $C_2^{\prime}$ for all irreducible representations.
Note also that we have omitted the eigenvalues for matrices in the irreducible representation $\Sigma_6^{(2)\prime}$. These eigenvalues are identical to those for $\Sigma_6^{(2)}$, except for elements in the conjugacy classes $C_8$, $C_8^{\prime}$, where for $\Sigma_6^{(2)\prime}$ the eigenvalues for elements in $C_8$, $C_8^{\prime}$ respectively are given by those for elements in $C_8^{\prime}$, $C_8$ respectively in the representation $\Sigma_6^{(2)}$.

From considering the set of eigenvalues $X_C$ for group elements in $C$ in the representation $\gamma_1^{(i)}$, $i=4,8$, we see that there is a choice of $(t_1^C,t_2^C) \in X_C$ such that $\mathcal{E}_{t_1^C,t_2^C} = X_C$, for all conjugacy classes $C$. We present one such choice of eigenvalues $(t_1^C,t_2^C)$ in Table \ref{table:t1,t2-IP3}.

\begin{figure}[tb]
\begin{center}
  \includegraphics[width=110mm]{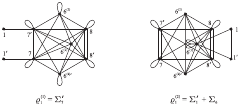}\\
 \caption{The McKay graphs $\mathcal{G}^{\varrho_1^{(i)}}_{PGL(2;7)}$, $i=1,2$.} \label{Fig-McKay_Graph-IP3-rho1}
\end{center}
\end{figure}

Thus we obtain two non-conjugate embeddings of $PGL(2;7)$ in $G_2$. We set $\varrho_1^{(1)} = \Sigma_7^{\prime}$, $\varrho_1^{(2)} = \Sigma_1^{\prime} + \Sigma_6^{(1)}$.
The McKay graphs $\mathcal{G}^{\varrho_1^{(1)}}_{PGL(2;7)}$, $\mathcal{G}^{\varrho_1^{(2)}}_{PGL(2;7)}$ for $\varrho_1^{(1)}$, $\varrho_1^{(2)}$ are given in Figure \ref{Fig-McKay_Graph-IP3-rho1}, where we use the same notation as previously.
Note that the graph given in \cite[Figure 2]{he:2003} is not the McKay graph for the restriction $\varrho_1$ of the fundamental seven-dimensional representation of $G_2$ as claimed in \cite{he:2003}, but is rather the McKay graph for the seven-dimensional representation $\Sigma_7$, which as shown above does not define an embedding of $PGL(2;7)$ in $G_2$.

\begin{table}[tb]
\begin{center}
\begin{tabular}{|c||c|c|c|c|c|c|c|c|} \hline
& \multicolumn{2}{|c|}{$\varrho_1^{(1)}$} & \multicolumn{2}{|c|}{$\varrho_1^{(2)}$} & & & & \\
$C$ & $(t_1^C,t_2^C)$ & $\chi_{\varrho_2^{(1)}}(C)$ & $(t_1^C,t_2^C)$ & $\chi_{\varrho_2^{(2)}}(C)$ & $\chi_{\gamma_2^{(1)}}(C)$ & $\chi_{\gamma_2^{(2)}}(C)$ & $\chi_{\gamma_2^{(3)}}(C)$ & $\chi_{\gamma_2^{(4)}}(C)$ \\
\hline \hline $C_1$ & $(1,1)$ & 14 & $(1,1)$ & 14 & 14 & 14 & 14 & 14 \\
\hline $C_2$ & $(1,-1)$ & -2 & $(1,-1)$ & -2 & -2 & -2 & 2 & 2 \\
\hline $C_2^{\prime}$ & $(1,-1)$ & -2 & $(1,-1)$ & -2 & 2 & -2 & 2 & -2 \\
\hline $C_8$ & $(i,\nu^7)$ & 0 & $(-1,\nu)$ & 0 & 0 & 0 & $\sqrt{2}$ & $\sqrt{2}$ \\
\hline $C_8^{\prime}$ & $(i,\nu^5)$ & 0 & $(-1,\nu)$ & 0 & 0 & 0 & $-\sqrt{2}$ & $-\sqrt{2}$ \\
\hline $C_4$ & $(-1,i)$ & 2 & $(1,i)$ & 2 & 2 & 2 & 0 & 0 \\
\hline $C_7$ & $(\eta,\eta^5)$ & 0 & $(\eta,\eta^5)$ & 0 & 0 & 0 & 0 & 0 \\
\hline $C_6$ & $(-1,\mu)$ & 1 & $(-1,\mu)$ & 1 & -1 & 1 & -1 & 1 \\
\hline $C_3$ & $(1,\omega)$ & -1 & $(1,\omega)$ & -1 & -1 & -1 & -1 & -1 \\
\hline
\end{tabular} \\
\caption{Choice of eigenvalues $(t_1^C,t_2^C)$ for $\varrho_1^{(i)}$, $i=1,2$, and corresponding values of $\chi_{\varrho_2^{(i)}}(C)$.} \label{table:t1,t2-IP3}
\end{center}
\end{table}

Since $\chi_{\varrho_2^{(i)}}(C) = \Phi_2(t_1^C,t_2^C)$, we see from Table \ref{table:t1,t2-II2} that the decomposition of the fundamental fourteen-dimensional representation into irreducible representations of $PGL(2;7)$ is given by $\varrho_2 = \gamma_2^{(2)} = \Sigma_6^{(1)} + \Sigma_8^{\prime}$ for both $i=1,2$.
The values of $x^{(i)} = \chi_{\varrho_1^{(i)}}(C) \in [-2,7]$, $y = \chi_{\varrho_2}(C) \in [-2,14]$ for $PGL(2;7)$ are given in Table \ref{table:(x,y)-IP3}.

\begin{table}[tb]
\begin{center}
\begin{tabular}{|c||c|c|c|c|c|c|c|c|c|} \hline
$C$ & $C_1$ & $C_2$ & $C_2^{\prime}$ & $C_8$ & $C_8^{\prime}$ & $C_4$ & $C_7$ & $C_6$ & $C_3$ \\
\hline $\chi_{\varrho_1^{(1)}}(C) \in [-2,7]$ & 7 & -1 & -1 & 1 & 1 & -1 & 0 & -1 & 1 \\
\hline $\chi_{\varrho_1^{(2)}}(C) \in [-2,7]$ & 7 & -1 & -1 & -1 & -1 & 3 & 0 & -1 & 1 \\
\hline $\chi_{\varrho_2}(C) \in [-2,14]$ & 14 & -2 & -2 & 0 & 0 & 2 & 0 & 1 & -1 \\
\hline
\end{tabular} \\
\caption{$\chi_{\varrho_j}(C)$ for group $PGL(2;7)$, $j=1,2$.} \label{table:(x,y)-IP3}
\end{center}
\end{table}

Then from (\ref{eqn:moments-subgroupG2}) and Tables \ref{table:Character_table-IP3}, \ref{table:(x,y)-IP3} and \ref{Table:subgroupsG2-orbits(theta1,theta2)}, we see that
\begin{align*}
\varsigma_{m,n} & = \frac{1}{336} \Omega^W(0,0) + \frac{21+28}{336} \Omega^W(0,1/2) + \frac{56}{336} \Omega^W(0,1/3) \\
& \quad + \frac{56}{336} \Omega^W(1/6,1/2) + \frac{48}{336} \Omega^W(1/7,3/7) + \frac{42}{336} \Omega' + \frac{42+42}{336} \Omega'',
\end{align*}
where $\Omega^W(\theta_1,\theta_2)$ is as in Section \ref{sect:II1} and $(\Omega',\Omega'')$ is $(\Omega^W(1/4,1/2),\Omega^W(1/8,3/8))$ for $\varrho_1^{(1)}$, and $(\Omega^W(0,1/4),\Omega^W(1/8,1/2))$ for $\varrho_1^{(2)}$.
Then we obtain:

\begin{Thm}
The joint spectral measure (over $\bbT^2$) for the non-conjugate embeddings of $PGL(2;7)$ into the fundamental representations of $G_2$ is
\begin{equation}
\begin{split}
\mathrm{d}\varepsilon & = \frac{1}{1344\pi^4} J^2 \, \mathrm{d}_7 \times \mathrm{d}_7 + \frac{1}{1152\pi^4} J^2 \, \mathrm{d}_6 \times \mathrm{d}_6 + \frac{1}{48} K' \, \mathrm{d}_4 \times \mathrm{d}_4 + \frac{1}{4} \, \mathrm{d}_3 \times \mathrm{d}_3 \\
& \quad + \frac{1}{36} \beta \, \mathrm{d}_2 \times \mathrm{d}_2 - \frac{23}{504} \, \mathrm{d}_1 \times \mathrm{d}_1 + \frac{1}{48} \sum_{g \in D_{12}} \delta_{g(\mu,t)},
\end{split}
\end{equation}
for $\mu$ a primitive 8-th root of unity, where for the embedding of $PGL(2;7)$ in $G_2$ given by $\varrho_1^{(1)} = \Sigma_7^{\prime}$ we have $K' = 16-K$, $\beta = 4$ and $t=\mu^3$, whilst for the embedding given by $\varrho_1^{(2)} = \Sigma_1^{\prime} + \Sigma_6$ we have $K' = K$, $\beta = 7$ and $t=-1$, with $K$ as in Theorem \ref{thm:measureII1},
and where $\mathrm{d}_m$ is the uniform measure over $m^{\mathrm{th}}$ roots of unity and $\delta_x$ is the Dirac measure at the point $x$.
\end{Thm}

\subsection{Group $PSL(2;8)$}

The subgroup $PSL(2;8)$ of $G_2$ is an irreducible primitive group of order 504.
It has nine irreducible representations, all real, five of which have dimension less than or equal to 7. The character table for $PSL(2;8)$ is given in
Table \ref{table:Character_table-IP2} \cite{james/liebeck:2001} (the orders of the elements in each conjugacy class can be obtained from \cite{lopez_pena/majid/rietsch:2010}).
The spectral measure for the first three seven-dimensional representations are equal, since the conjugacy classes $C_{9}$, $C_{9}^{\prime}$, $C_{9}^{\prime\prime}$ each have the same order.

\begin{table}[tb]
\begin{center}
\begin{tabular}{|c||c|c|c|c|c|c|c|c|c|} \hline
$C$ & $C_1$ & $C_3$ & $C_9$ & $C_9^{\prime}$ & $C_9^{\prime\prime}$ & $C_2$ & $C_7$ & $C_7^{\prime}$ & $C_7^{\prime\prime}$ \\
\hline $|C|$ & 1 & 56 & 56 & 56 & 56 & 63 & 72 & 72 & 72 \\
\hline \hline $\Sigma_1$ & 1 & 1 & 1 & 1 & 1 & 1 & 1 & 1 & 1 \\
\hline $\Sigma_7^{(1)}$ & 7 & 1 & $-p$ & $-q$ & $p+q$ & -1 & 0 & 0 & 0 \\
\hline $\Sigma_7^{(1)\prime}$ & 7 & 1 & $p+q$ & $-p$ & $-q$ & -1 & 0 & 0 & 0 \\
\hline $\Sigma_7^{(1)\prime\prime}$ & 7 & 1 & $-q$ & $p+q$ & $-p$ & -1 & 0 & 0 & 0 \\
\hline $\Sigma_7^{(2)}$ & 7 & -2 & 1 & 1 & 1 & -1 & 0 & 0 & 0 \\
\hline $\Sigma_8$ & 8 & -1 & -1 & -1 & -1 & 0 & 1 & 1 & 1 \\
\hline $\Sigma_9$ & 9 & 0 & 0 & 0 & 0 & 1 & $2\cos(2\pi/7)$ & $2\cos(4\pi/7)$ & $2\cos(6\pi/7)$ \\
\hline $\Sigma_9^{\prime}$ & 9 & 0 & 0 & 0 & 0 & 1 & $2\cos(4\pi/7)$ & $2\cos(6\pi/7)$ & $2\cos(2\pi/7)$ \\
\hline $\Sigma_9^{\prime\prime}$ & 9 & 0 & 0 & 0 & 0 & 1 & $2\cos(6\pi/7)$ & $2\cos(2\pi/7)$ & $2\cos(4\pi/7)$ \\
\hline
\end{tabular} \\
\caption{Character table for $PSL(2;8)$, where $p=2\cos(4\pi/9)$, $q=2\cos(8\pi/9)$.} \label{table:Character_table-IP2}
\end{center}
\end{table}

We thus have four candidates for the restriction $\varrho_1$ of the fundamental representation $\rho_1$ of $G_2$ to $PSL(2;8)$. The Kronecker squares of $\Sigma_7^{(1)}$, $\Sigma_7^{(2)}$ decompose into irreducibles as
$$(\Sigma_7^{(1)})^2 = \mathrm{id} + \Sigma_7^{(1)} + \Sigma_7^{(1)\prime} + \Sigma_7^{(2)} + \Sigma_9 + \Sigma_9^{\prime} + \Sigma_9^{\prime\prime}, \qquad
(\Sigma_7^{(2)})^2 = \mathrm{id} + \Sigma_7^{(1)} + \Sigma_7^{(1)\prime} + \Sigma_7^{(1)\prime\prime} + \Sigma_9 + \Sigma_9^{\prime} + \Sigma_9^{\prime\prime},$$
where $\mathrm{id} = \Sigma_1$. The irreducible representation $\Sigma_7^{(2)}$ does not appear in the decomposition of its Kronecker square into irreducibles, and therefore does not give an embedding of $PSL(2;8)$ into $G_2$.

\begin{figure}[tb]
\begin{center}
  \includegraphics[width=55mm]{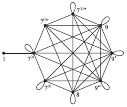}\\
 \caption{The McKay graph $\mathcal{G}^{\varrho_1^{(1)}}_{PSL(2;8)}$.} \label{Fig-McKay_Graph-IP2-rho1}
\end{center}
\end{figure}

Thus we obtain that there are precisely three seven-dimensional representations give non-conjugate embeddings of $PSL(2;8)$ into $G_2$. Let $\varrho_1^{(1)} = \Sigma_7^{(1)}$, $\varrho_1^{(2)} = \Sigma_7^{(1)\prime}$ and $\varrho_1^{(3)} = \Sigma_7^{(1)\prime\prime}$.
The McKay graph for $\varrho_1^{(1)}$ is given in Figure \ref{Fig-McKay_Graph-IP2-rho1}, where we use the same notation as previously. The McKay graphs for $\varrho_1^{(j)}$, $j=2,3$, are similar.
Note that the graph given in \cite[Figure 2]{he:2003} is not the McKay graph for a restriction $\varrho_1^{(j)}$, $j \in \{1,2,3\}$, of the fundamental seven-dimensional representation of $G_2$ as claimed in \cite{he:2003}, but is rather the McKay graph for the seven-dimensional representation $\Sigma_7^{(2)}$.

From the decomposition of the Kronecker square of $\varrho_1^{(j)}$ and dimension considerations there is only one possibility for the fourteen-dimensional representation $\varrho_2$, that is, $\varrho_2 = \Sigma_7^{(1)\prime} + \Sigma_7^{(2)}$.
The values of $x = \chi_{\varrho_1^{(j)}}(C) \in [-2,7]$, $j \in \{1,2,3\}$, and $y = \chi_{\varrho_2}(C) \in [-2,14]$ for $PSL(2;8)$ are given in Table \ref{table:(x,y)-IP2}.

\begin{table}[tb]
\begin{center}
\begin{tabular}{|c||c|c|c|c|c|c|c|c|c|} \hline
$C$ & $C_1$ & $C_3$ & $C_9$ & $C_9^{\prime}$ & $C_9^{\prime\prime}$ & $C_2$ & $C_7$ & $C_7^{\prime}$ & $C_7^{\prime\prime}$ \\
\hline $\chi_{\varrho_1^{(j)}}(C) \in [-2,7]$ & 7 & 1 & $-p$ & $-q$ & $p+q$ & -1 & 0 & 0 & 0 \\
\hline $\chi_{\varrho_2}(C) \in [-2,14]$ & 14 & -1 & $p+q+1$ & $1-p$ & $1-q$ & -2 & 0 & 0 & 0 \\
\hline
\end{tabular} \\
\caption{$\chi_{\varrho_1^{(j)}}(C)$ and $\chi_{\varrho_2}(C)$ for group $PSL(2;8)$, $j \in \{1,2,3\}$.} \label{table:(x,y)-IP2}
\end{center}
\end{table}

Then from (\ref{eqn:moments-subgroupG2}) and Tables \ref{table:Character_table-IP2} and \ref{Table:subgroupsG2-orbits(theta1,theta2)}, we see that
\begin{align*}
\varsigma_{m,n} & = \frac{1}{504} \Omega^W(0,0) + \frac{56}{504} \Omega^W(0,1/3) + \frac{63}{504} \Omega^W(0,1/2) + \frac{72+72+72}{504} \Omega^W(1/7,3/7) \\
& \quad + \frac{56}{504} \Omega^W(1/9,4/9) + \frac{56}{504} \Omega^W(1/9,1/3) + \frac{56}{504} \Omega^W(2/9,5/9),
\end{align*}
where $\Omega^W(\theta_1,\theta_2)$ is as in Section \ref{sect:II1}.
Then we obtain:

\begin{Thm}
The joint spectral measure (over $\bbT^2$) for all embeddings of $PSL(2;8)$ into the fundamental representations of $G_2$ is
\begin{equation}
\begin{split}
\mathrm{d}\varepsilon & = \frac{1}{448\pi^4} J^2 \, \mathrm{d}_7 \times \mathrm{d}_7 + \frac{1}{6} \, \mathrm{d}_3 \times \mathrm{d}_3 + \frac{1}{6} \, \mathrm{d}_2 \times \mathrm{d}_2 - \frac{5}{126} \, \mathrm{d}_1 \times \mathrm{d}_1 - \frac{1}{18} \, \mathrm{d}^{(1)} \\
& \quad + \frac{1}{864\pi^4} a_3^{-1} J^2 \, \mathrm{d}^{((9))} + \frac{1}{864\pi^4} a_1^{-1} J^2 \, \mathrm{d}^{((9/2))} + \frac{1}{864\pi^4} a_2^{-1} J^2 \, \mathrm{d}^{((9/4))},
\end{split}
\end{equation}
where $\mathrm{d}^{((n))}$ is as in Definition \ref{def:4measures}, $\mathrm{d}_m$ is the uniform measure over $m^{\mathrm{th}}$ roots of unity and $\mathrm{d}^{(k+4)}$ is the uniform measure on the points in $C_k^W$.
\end{Thm}

\subsection{Group $PSL(2;13)$}

The subgroup $PSL(2;13)$ of $G_2$ is an irreducible primitive group of order 1092. It has nine irreducible representations, all real, and its character table is given in Table \ref{table:Character_table-IP1} \cite{cohen:1998}.

\begin{table}[tb]
\begin{center}
\begin{tabular}{|c||c|c|c|c|c|c|c|c|c|} \hline
$C$ & $C_1$ & $C_{13}$ & $C_{13}^{\prime}$ & $C_2$ & $C_7$ & $C_7^{\prime}$ & $C_7^{\prime\prime}$ & $C_6$ & $C_3$ \\
\hline $|C|$ & 1 & 84 & 84 & 91 & 156 & 156 & 156 & 182 & 182 \\
\hline \hline $\Sigma_1$ & 1 & 1 & 1 & 1 & 1 & 1 & 1 & 1 & 1 \\
\hline $\Sigma_7$ & 7 & $p_-$ & $p_+$ & -1 & 0 & 0 & 0 & -1 & 1 \\
\hline $\Sigma_7^{\prime}$ & 7 & $p_+$ & $p_-$ & -1 & 0 & 0 & 0 & -1 & 1 \\
\hline $\Sigma_{12}$ & 12 & -1 & -1 & 0 & $r_1$ & $r_2$ & $r_3$ & 0 & 0 \\
\hline $\Sigma_{12}^{\prime}$ & 12 & -1 & -1 & 0 & $r_2$ & $r_3$ & $r_1$ & 0 & 0 \\
\hline $\Sigma_{12}^{\prime\prime}$ & 12 & -1 & -1 & 0 & $r_3$ & $r_1$ & $r_2$ & 0 & 0 \\
\hline $\Sigma_{13}$ & 13 & 0 & 0 & 1 & -1 & -1 & -1 & 1 & 1 \\
\hline $\Sigma_{14}$ & 14 & 1 & 1 & -2 & 0 & 0 & 0 & 1 & -1 \\
\hline $\Sigma_{14}^{\prime}$ & 14 & 1 & 1 & 2 & 0 & 0 & 0 & -1 & -1 \\
\hline
\end{tabular} \\
\caption{Character table for $PSL(2;13)$, where $p_{\pm} = (1\pm\sqrt{13})/2$, $r_j=-2\cos(2j\pi/7)$.} \label{table:Character_table-IP1}
\end{center}
\end{table}

Only three of the irreducible representation have dimension less than or equal to 7. These are the identity representation, and two seven-dimensional representations $\Sigma_7$, $\Sigma_7^{\prime}$ whose character values only differ on the two conjugacy classes $C_{13}$, $C_{13}^{\prime}$ whose elements have order 13. Here $\chi_{\Sigma_7}(g) = p,q$ for $g \in C_{13},C_{13}^{\prime}$ respectively, whilst $\chi_{\Sigma_7^{\prime}}(g) = q,p$ respectively, where $p=(1+\sqrt{13})/2=1+\zeta+\zeta^3+\zeta^4+\zeta^9+\zeta^{10}+\zeta^{12}$, $q=(1-\sqrt{13})/2=1+\zeta^2+\zeta^5+\zeta^6+\zeta^7+\zeta^8+\zeta^{11}$, for $\zeta = e^{2\pi i/13}$.
Thus the spectral measure for both seven-dimensional representations are equal, since $C_{13}$, $C_{13}^{\prime}$ both have the same order.

\begin{figure}[tb]
\begin{center}
  \includegraphics[width=60mm]{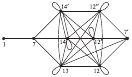} \hspace{15mm} \includegraphics[width=60mm]{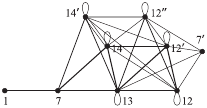}\\
 \caption{Two presentations of the McKay graph $\mathcal{G}^{\varrho_1^{(j)}}_{PSL(2;13)}$, $j \in \{ 1,2 \}$.} \label{Fig-McKay_Graph-IP1-rho1}
\end{center}
\end{figure}

There are thus two non-conjugate embeddings $\varrho_1^{(1)} := \Sigma_7$ and $\varrho_1^{(2)} := \Sigma_7^{\prime}$ of $PSL(2;13)$ in $G_2$.
The McKay graph for the fundamental seven-dimensional representation $\varrho_1^{(j)}$, $j \in \{ 1,2 \}$,  is given in \cite[Figure 2]{he:2003}, and is reproduced (twice) in Figure \ref{Fig-McKay_Graph-IP1-rho1} for completeness, where we use the same notation as previously. The figure on the right hand side illustrates the resemblance of $\mathcal{G}^{\varrho_1^{(j)}}_{PSL(2;13)}$ with the McKay graph of $G_2$ itself.

The set $X_C$ of eigenvalues of group elements in each conjugacy class for $\varrho_1^{(1)} = \Sigma_7$ are given in Table \ref{table:evalues-IP1}, along with a choice of $(t_1^C,t_2^C) \in X_C$ such that $\mathcal{E}_{t_1^C,t_2^C} = X_C$. The eigenvalues for $\varrho_1^{(2)}$ are similar.
The decomposition of the Kronecker square of $\varrho_1^{(j)}$, $j \in \{1,2\}$, into irreducibles is given by
$$\left(\varrho_1^{(j)}\right)^2 = \mathrm{id} + \varrho_1^{(j)} + \Sigma_{13} + \Sigma_{14} + \Sigma_{14}^{\prime},$$
where as before the notation $\Sigma_n$ denotes an irreducible representation of $PSL(2;13)$ of dimension $n$.
From dimension considerations there are thus two candidates, $\Sigma_{14}$ and $\Sigma_{14}^{\prime}$, for the fourteen-dimensional representation $\varrho_2$, which are given by the two fourteen-dimensional irreducible representations.
However, since $\chi_{\varrho_2}(C) = \Phi_2(t_1^C,t_2^C)$, we see from Table \ref{table:evalues-IP1} that the decomposition of the fundamental fourteen-dimensional representation into irreducible representations of $PSL(2;13)$ is given by $\varrho_2 = \Sigma_{14}$, not $\Sigma_{14}^{\prime}$.

\begin{table}[tb]
\begin{center}
\begin{tabular}{|c||c|c|c|c|c|c|} \hline
$C$ & $X_C$ & $(t_1^C,t_2^C)$ & $\chi_{\varrho_1^{(1)}}(C)$ & $\chi_{\varrho_2}(C)$ & $\chi_{\Sigma_{14}}(C)$ & $\chi_{\Sigma_{14}^{\prime}}(C)$ \\
\hline \hline $C_1$ & $(1,1,1,1,1,1,1)$ & $(1,1)$ & 7 & 14 & 14 & 14 \\
\hline $C_{13}$ & $(1,\zeta,\zeta^3,\zeta^4,\zeta^9,\zeta^{10},\zeta^{12})$ & $(\zeta,\zeta^4)$ & $(1+\sqrt{13})/2$ & 1 & 1 & 1 \\
\hline $C_{13}^{\prime}$ & $(1,\zeta^2,\zeta^5,\zeta^6,\zeta^7,\zeta^8,\zeta^{11})$ & $(\zeta^2,\zeta^7)$ & $(1-\sqrt{13})/2$ & 1 & 1 & 1 \\
\hline $C_2$ & $(1,1,1,-1,-1,-1,-1)$ & $(1,-1)$ & -1 & -2 & -2 & 2\\
\hline $C_7$ & $(1,\eta,\eta^2,\eta^3,\eta^4,\eta^5,\eta^6)$ & $(\eta,\eta^5)$ & 0 & 0 & 0 & 0 \\
\hline $C_7^{\prime}$ & $(1,\eta,\eta^2,\eta^3,\eta^4,\eta^5,\eta^6)$ & $(\eta,\eta^5)$ & 0 & 0 & 0 & 0 \\
\hline $C_7^{\prime\prime}$ & $(1,\eta,\eta^2,\eta^3,\eta^4,\eta^5,\eta^6)$ & $(\eta,\eta^5)$ & 0 & 0 & 0 & 0 \\
\hline $C_6$ & $(1,-1,-1,\mu,\mu^2,\mu^4,\mu^5)$ & $(-1,\mu)$ & -1 & 1 & 1 & -1 \\
\hline $C_3$ & $(1,1,1,\omega,\omega,\omega^2,\omega^2)$ & $(1,\omega)$ & 1 & -1 & -1 & -1\\
\hline
\end{tabular} \\
\caption{Choice of eigenvalues $(t_1^C,t_2^C)$ for $\varrho_1^{(1)}$ and corresponding values of $\chi_{\varrho_2}(C)$, where $\omega = e^{2\pi i/3}$, $\mu = e^{2\pi i/6}$, $\eta = e^{2\pi i/7}$ and $\zeta = e^{2\pi i/13}$.} \label{table:evalues-IP1}
\end{center}
\end{table}

Then from (\ref{eqn:moments-subgroupG2}) and Tables \ref{table:evalues-IP1} and \ref{Table:subgroupsG2-orbits(theta1,theta2)}, we see that
\begin{align*}
\varsigma_{m,n} & = \frac{1}{1092} \Omega^W(0,0) + \frac{91}{1092} \Omega^W(0,1/2) + \frac{182}{1092} \Omega^W(0,1/3) + \frac{182}{1092} \Omega^W(1/6,1/2) \\
& \quad + \frac{156+156+156}{1092} \Omega^W(1/7,3/7) + \frac{84}{1092} \Omega^W(1/13,4/13)  + \frac{84}{1092} \Omega^W(2/13,7/13),
\end{align*}
where $\Omega^W(\theta_1,\theta_2)$ is as in Section \ref{sect:II1}.
Then we obtain:

\begin{Thm}
The joint spectral measure (over $\bbT^2$) for all embeddings of $PSL(2;13)$ into the fundamental representations of $G_2$ is
\begin{equation}
\begin{split}
\mathrm{d}\varepsilon & = \frac{1}{448\pi^4} J^2 \, \mathrm{d}_7 \times \mathrm{d}_7 + \frac{1}{1152\pi^4} J^2 \, \mathrm{d}_6 \times \mathrm{d}_6 + \frac{1}{4} \, \mathrm{d}_3 \times \mathrm{d}_3 + \frac{1}{9} \, \mathrm{d}_2 \times \mathrm{d}_2 \\
& \quad - \frac{22}{819} \, \mathrm{d}_1 \times \mathrm{d}_1 - \frac{1}{12} \, \mathrm{d}^{(1)} + \frac{1}{156} \sum_{g \in D_{12}} (\delta_{g(\zeta,\zeta^4)} + \delta_{g(\zeta^2,\zeta^7)})
\end{split}
\end{equation}
for $\zeta$ a primitive 13-th root of unity, where $\mathrm{d}_m$ is the uniform measure over $m^{\mathrm{th}}$ roots of unity, $\mathrm{d}^{(k+4)}$ is the uniform measure on the points in $C_k^W$, and $\delta_x$ is the Dirac measure at the point $x$.
\end{Thm}

\subsection{Group $PU(3;3) \cong G_2(2)'$}

The subgroup $PU(3;3)$ of $G_2$ is an irreducible primitive group of order 6048.
It has fourteen irreducible representations (seven real repreentations, one quaternionic representation $\Sigma_6$, and three pairs of complex conjugate representations) and its character table is given in Table \ref{table:Character_table-IP4} \cite{conway/curtis/norton/parker/wilson:1985}.

\begin{table}[tb]
\begin{center}
\begin{tabular}{|@{\hspace{1.5mm}}c@{\hspace{1.5mm}}||@{\hspace{1mm}}c@{\hspace{1mm}}|@{\hspace{1mm}}c@{\hspace{1mm}}|@{\hspace{1mm}}c @{\hspace{1mm}}|@{\hspace{1mm}}c@{\hspace{1mm}}|@{\hspace{1mm}}c@{\hspace{1mm}}|@{\hspace{1mm}}c@{\hspace{1mm}}|@{\hspace{1mm}}c @{\hspace{1mm}}|@{\hspace{1mm}}c@{\hspace{1mm}}|@{\hspace{1mm}}c@{\hspace{1mm}}|@{\hspace{1mm}}c@{\hspace{1mm}}|@{\hspace{1mm}}c @{\hspace{1mm}}|@{\hspace{1mm}}c@{\hspace{1mm}}|@{\hspace{1mm}}c@{\hspace{1mm}}|@{\hspace{1mm}}c@{\hspace{1mm}}|} \hline
$C$ & $C_1$ & $C_3$ & $C_2$ & $C_4$ & $C_4^{\prime}$ & $C_4^{\prime\prime}$ & $C_{12}$ & $C_{12}^{\prime}$ & $C_6$ & $C_3^{\prime}$ & $C_8$ & $C_8^{\prime}$ & $C_7$ & $C_7^{\prime}$ \\
\hline $|C|$ & 1 & 56 & 63 & 63 & 63 & 378 & 504 & 504 & 504 & 672 & 756 & 756 & 864 & 864 \\
\hline \hline $\Sigma_1$ & 1 & 1 & 1 & 1 & 1 & 1 & 1 & 1 & 1 & 1 & 1 & 1 & 1 & 1 \\
\hline $\Sigma_6$ & 6 & -3 & -2 & -2 & -2 & 2 & 1 & 1 & 1 & 0 & 0 & 0 & -1 & -1 \\
\hline $\Sigma_7$ & 7 & -2 & 3 & $-1-2i$ & $-1+2i$ & 1 & $-1-i$ & $-1+i$ & 0 & 1 & $-i$ & $i$ & 0 & 0 \\
\hline $\Sigma_7^{\ast}$ & 7 & -2 & 3 & $-1+2i$ & $-1-2i$ & 1 & $-1+i$ & $-1-i$ & 0 & 1 & $i$ & $-i$ & 0 & 0 \\
\hline $\Sigma_7^{\prime}$ & 7 & -2 & -1 & 3 & 3 & -1 & 0 & 0 & 2 & 1 & -1 & -1 & 0 & 0 \\
\hline $\Sigma_{14}$ & 14 & 5 & -2 & 2 & 2 & 2 & -1 & -1 & 1 & -1 & 0 & 0 & 0 & 0 \\
\hline $\Sigma_{21}$ & 21 & 3 & 1 & $-3-2i$ & $-3+2i$ & -1 & $-i$ & $i$ & 1 & 0 & $i$ & $-i$ & 0 & 0 \\
\hline $\Sigma_{21}^{\ast}$ & 21 & 3 & 1 & $-3+2i$ & $-3-2i$ & -1 & $i$ & $-i$ & 1 & 0 & $-i$ & $i$ & 0 & 0 \\
\hline $\Sigma_{21}^{\prime}$ & 21 & 3 & 5 & 1 & 1 & 1 & 1 & 1 & -1 & 0 & -1 & -1 & 0 & 0 \\
\hline $\Sigma_{27}$ & 27 & 0 & 3 & 3 & 3 & -1 & 0 & 0 & 0 & 0 & 1 & 1 & -1 & -1 \\
\hline $\Sigma_{28}$ & 28 & 1 & -4 & $-4i$ & $4i$ & 0 & $i$ & $-i$ & -1 & 1 & 0 & 0 & 0 & 0 \\
\hline $\Sigma_{28}^{\ast}$ & 28 & 1 & -4 & $4i$ & $-4i$ & 0 & $-i$ & $i$ & -1 & 1 & 0 & 0 & 0 & 0 \\
\hline $\Sigma_{32}$ & 32 & -4 & 0 & 0 & 0 & 0 & 0 & 0 & 0 & -1 & 0 & 0 & $\frac{1+i\sqrt{7}}{2}$ & $\frac{1-i\sqrt{7}}{2}$ \\
\hline $\Sigma_{32}^{\ast}$ & 32 & -4 & 0 & 0 & 0 & 0 & 0 & 0 & 0 & -1 & 0 & 0 & $\frac{1-i\sqrt{7}}{2}$ & $\frac{1+i\sqrt{7}}{2}$ \\
\hline
\end{tabular} \\
\caption{Character table for $PU(3;3)$.} \label{table:Character_table-IP4}
\end{center}
\end{table}

There are two non-conjugate real seven-dimensional representations, $\gamma_1^{(1)} = \Sigma_1 + \Sigma_6$ and $\gamma_1^{(2)} = \Sigma_7^{\prime}$.
These both satisfy the condition that $\gamma_1^{(i)}$ appears in the decomposition of $(\gamma_1^{(i)})^2$.
We thus consider the eigenvalues of the representation matrices to determine which of the $\gamma_1^{(i)}$ are embeddings of $PU(3;3)$ in $G_2$.
These eigenvalues are given in Table \ref{table:evalues-IP4} for the representations $\Sigma_1$, $\Sigma_6$ and $\Sigma_7^{\prime}$. As described in Section \ref{sect:subgroupsG2}, these eigenvalues can be determined from the character table of $PU(3;3)$. The additional information that is needed is to note that the eigenvalues for group elements in $C_4$, $C_4^{\prime}$ and $C_4^{\prime\prime}$ all square to those for elements in $C_2$, those for $C_8$, $C_8^{\prime}$ square to those for $C_4$, $C_4^{\prime}$ respectively, those for $C_6$ square to those for $C_3$ and also cube to those for $C_2$, those for $C_{12}$ square to those for $C_6$ and cube to those for $C_4$ whilst those for $C_{12}^{\prime}$ also square to those for $C_6$ but cube to those for $C_4^{\prime}$ (see \cite{conway/curtis/norton/parker/wilson:1985}).

\begin{table}[tb]
\begin{center}
\begin{tabular}{|c||c|c|c||c|} \hline
 & $\Sigma_1$ & $\Sigma_6$ & $\Sigma_7^{\prime}$ & $(t_1^C,t_2^C)$ \\
\hline \hline $C_1$ & 1 & $(1,1,1,1,1,1)$ & $(1,1,1,1,1,1,1)$ & $(1,1)$ \\
\hline $C_3$ & 1 & $(\omega,\omega,\omega,\omega^2,\omega^2,\omega^2)$ & $(1,\omega,\omega,\omega,\omega^2,\omega^2,\omega^2)$  & $(\omega,\omega^2)$ \\
\hline $C_2$ & 1 & $(1,1,-1,-1,-1,-1)$ & $(1,1,1,-1,-1,-1,-1)$  & $(1,-1)$ \\
\hline $C_4$ & 1 & $(-1,-1,i,i,-i,-i)$ & $(1,1,1,i,i,-i,-i)$  & $(1,i)$ \\
\hline $C_4^{\prime}$ & 1 & $(-1,-1,i,i,-i,-i)$ & $(1,1,1,i,i,-i,-i)$  & $(1,i)$ \\
\hline $C_4^{\prime\prime}$ & 1 & $(1,1,i,i,-i,-i)$ & $(1,-1,-1,i,i,-i,-i)$  & $(-1,i)$ \\
\hline $C_{12}$ & 1 & $(\xi,\xi^2,\xi^5,\xi^7,\xi^{10},\xi^{11})$ & $(1,\xi,\xi^4,\xi^5,\xi^7,\xi^8,\xi^{11})$  & $(\xi,\xi^5)$ \\
\hline $C_{12}^{\prime}$ & 1 & $(\xi,\xi^2,\xi^5,\xi^7,\xi^{10},\xi^{11})$ & $(1,\xi,\xi^4,\xi^5,\xi^7,\xi^8,\xi^{11})$  & $(\xi,\xi^5)$ \\
\hline $C_6$ & 1 & $(\mu,\mu,\mu^2,\mu^4,\mu^5,\mu^5)$ & $(1,\mu,\mu,\mu^2,\mu^4,\mu^5,\mu^5)$  & $(\mu,\mu^2)$ \\
\hline $C_3^{\prime}$ & 1 & $(1,1,\omega,\omega,\omega^2,\omega^2)$ & $(1,1,1,\omega,\omega,\omega^2,\omega^2)$  & $(1,\omega)$ \\
\hline $C_8$ & 1 & $(i,-i,\nu,\nu^3,\nu^5,\nu^7)$ & $(1,-1,-1,\nu,\nu^3,\nu^5,\nu^7)$  & $(-1,\nu)$ \\
\hline $C_8^{\prime}$ & 1 & $(i,-i,\nu,\nu^3,\nu^5,\nu^7)$ & $(1,-1,-1,\nu,\nu^3,\nu^5,\nu^7)$  & $(-1,\nu)$ \\
\hline $C_7$ & 1 & $(\eta,\eta^2,\eta^3,\eta^4,\eta^5,\eta^6)$ & $(1,\eta,\eta^2,\eta^3,\eta^4,\eta^5,\eta^6)$  & $(\eta,\eta^5)$ \\
\hline $C_7^{\prime}$ & 1 & $(\eta,\eta^2,\eta^3,\eta^4,\eta^5,\eta^6)$ & $(1,\eta,\eta^2,\eta^3,\eta^4,\eta^5,\eta^6)$  & $(\eta,\eta^5)$ \\
\hline
\end{tabular} \\
\caption{Eigenvalues of group elements in each conjugacy class of $PU(3;3)$ for the irreducible representations $\Sigma_1$, $\Sigma_6$ and $\Sigma_7^{\prime}$, where $\omega = e^{2\pi i/3}$, $\mu = e^{2\pi i/6}$, $\eta = e^{2\pi i/7}$, $\nu = e^{2\pi i/8}$ and $\xi = e^{2\pi i/12}$} \label{table:evalues-IP4}
\end{center}
\end{table}

From considering the set of eigenvalues $X_C$ for group elements in $C$ in the representation $\gamma_1^{(i)}$, we see that there is no choice of $(t_1^C,t_2^C) \in X_C$ such that $\mathcal{E}_{t_1^C,t_2^C} = X_C$ for $i=1$ when $C=C_{12},C_{12}^{\prime}$. However, such a choice does exist for all conjugacy classes $C$ for $i=2$, thus we have $\varrho_1 = \gamma_1^{(2)}$. We present one such choice of eigenvalues $(t_1^C,t_2^C)$ in the final column of Table \ref{table:evalues-IP4}.

\begin{figure}[tb]
\begin{center}
  \includegraphics[width=70mm]{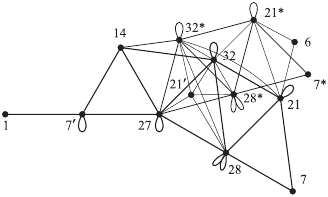}\\
 \caption{The McKay graphs $\mathcal{G}^{\varrho_1}_{PU(3;3)}$.} \label{Fig-McKay_Graph-IP4-rho1}
\end{center}
\end{figure}

Thus, up to conjugacy, there is only one embedding of $PU(3;3)$ in $G_2$.
The McKay graph $\mathcal{G}^{\varrho_1}_{PU(3;3)}$ is given in \cite[Figure 2]{he:2003}, and we reproduce it here in Figure \ref{Fig-McKay_Graph-IP4-rho1} for completeness.

The decomposition of the Kronecker square of $\varrho_1$ into irreducibles is given by
$$\varrho_1^2 = \mathrm{id} + \varrho_1 + \Sigma_{14} + \Sigma_{27},$$
where $\mathrm{id} = \Sigma_1$.
Thus the fourteen-dimensional representation $\varrho_2$ must be $\Sigma_{14}$, and we note that $\chi_{\varrho_2}(C) = \Phi_2(t_1^C,t_2^C)$ as required.

Then from (\ref{eqn:moments-subgroupG2}) and Tables \ref{table:Character_table-IP4} and \ref{Table:subgroupsG2-orbits(theta1,theta2)}, we see that
\begin{align*}
\varsigma_{m,n} & = \frac{1}{6048} \Omega^W(0,0) + \frac{56}{6048} \Omega^W(1/3,2/3) + \frac{63}{6048} \Omega^W(0,1/2) + \frac{672}{6048} \Omega^W(0,1/3) \\
& \quad + \frac{63+63}{6048} \Omega^W(0,1/4) + \frac{378}{6048} \Omega^W(1/4,1/2) + \frac{504}{6048} \Omega^W(1/6,1/3) \\
& \quad + \frac{864+864}{6048} \Omega^W(1/7,3/7) + \frac{756+756}{6048} \Omega^W(1/8,1/2) + \frac{504+504}{6048} \Omega^W(1/12,5/12),
\end{align*}
where $\Omega^W(\theta_1,\theta_2)$ is as in Section \ref{sect:II1}.
Then we obtain:

\begin{Thm}
The joint spectral measure (over $\bbT^2$) for all embeddings of $PU(3;3)$ into the fundamental representations of $G_2$ is
\begin{equation}
\begin{split}
\mathrm{d}\varepsilon & = \frac{1}{672\pi^4} J^2 \, \mathrm{d}_7 \times \mathrm{d}_7 + \frac{1}{144} (24-K) \, \mathrm{d}_4 \times \mathrm{d}_4 + \frac{1}{6} \, \mathrm{d}_3 \times \mathrm{d}_3 - \frac{1}{12} \, \mathrm{d}_2 \times \mathrm{d}_2 \\
& \quad + \frac{1}{168} \, \mathrm{d}_1 \times \mathrm{d}_1 + \frac{1}{1152\pi^4} J^2 \, \mathrm{d}^{(4)} + \frac{1}{6} \, \mathrm{d}^{(2)} - \frac{1}{12} \, \mathrm{d}^{(1)} + \frac{1}{48} \sum_{g \in D_{12}} \delta_{g(\mu,-1)},
\end{split}
\end{equation}
for $\mu$ \ primitive 8-th root of unity, where $K(\theta_1,\theta_2) = (\sin(2\pi(\theta_1+\theta_2))-\sin(2\pi(2\theta_1-\theta_2))-\sin(2\pi(2\theta_2-\theta_1)))^2$, $\mathrm{d}_m$ is the uniform measure over $m^{\mathrm{th}}$ roots of unity, $\mathrm{d}^{(k+4)}$ is the uniform measure on the points in $C_k^W$, and $\delta_x$ is the Dirac measure at the point $x$.
\end{Thm}

\subsection{Group $G_2(2)$} \label{sect:IP5}

The subgroup $G_2(2) = G_2(\mathbb{F}_2)$ of $G_2 = G_2(\mathbb{C})$ is an irreducible primitive group of order 12096. It is the group $G_2$ defined over the Galois field $\mathbb{F}_2$.
It has sixteen irreducible representations (fourteen real and two complex conjugate representations), and its character table is given in
\cite{he:2003} (the orders of the elements in each conjugacy class can be obtained from \cite{koca/koc:1994}).

\begin{table}[tb]
\begin{center}
\begin{tabular}{|@{\hspace{1.5mm}}c@{\hspace{1.5mm}}||@{\hspace{1mm}}c@{\hspace{1mm}}|@{\hspace{1mm}}c@{\hspace{1mm}}|@{\hspace{1mm}}c @{\hspace{1mm}}|@{\hspace{1mm}}c@{\hspace{1mm}}|@{\hspace{1mm}}c@{\hspace{1mm}}|@{\hspace{1mm}}c@{\hspace{1mm}}|@{\hspace{1mm}}c @{\hspace{1mm}}|@{\hspace{1mm}}c@{\hspace{1mm}}|@{\hspace{1mm}}c@{\hspace{1mm}}|@{\hspace{1mm}}c@{\hspace{1mm}}|@{\hspace{1mm}}c @{\hspace{1mm}}|@{\hspace{1mm}}c@{\hspace{1mm}}|@{\hspace{1mm}}c@{\hspace{1mm}}|@{\hspace{1mm}}c@{\hspace{1mm}}|@{\hspace{1mm}}c @{\hspace{1mm}}|@{\hspace{1mm}}c@{\hspace{1mm}}|} \hline
$C$ & $C_1$ & $C_3$ & $C_2$ & $C_4$ & $C_2^{\prime}$ & $C_4^{\prime}$ & $C_4^{\prime\prime}$ & $C_6$ & $C_3^{\prime}$ & $C_{12}$ & $C_{12}^{\prime}$ & $C_{12}^{\prime\prime}$ & $C_8$ & $C_8^{\prime}$ & $C_7$ & $C_6^{\prime}$ \\
\hline $|C|$                & 1  & 56 & 63 & 126 & 252 & 252 & 378 & 504 & 672 & 1008 & 1008 & 1008 & 1512 & 1512 & 1728 & 2016 \\
\hline $\Sigma_1$           & 1  & 1  & 1  & 1  & 1  & 1  & 1  & 1  & 1  & 1 & 1 & 1 & 1 & 1 & 1 & 1 \\
\hline $\Sigma_1'$          & 1  & 1  & 1  & 1  & -1 & -1 & 1  & 1  & 1  & -1 & -1 & 1 & -1 & 1 & 1 & -1 \\
\hline $\Sigma_6$           & 6  & -3 & -2 & -2 & 0  & 0  & 2  & 1  & 0  & $i\sqrt{3}$ & $-i\sqrt{3}$ & 1 & 0 & 0 & -1 & 0 \\
\hline $\Sigma_6^{\ast}$    & 6  & -3 & -2 & -2 & 0  & 0  & 2  & 1  & 0  & $-i\sqrt{3}$ & $i\sqrt{3}$ & 1 & 0 & 0 & -1 & 0 \\
\hline $\Sigma_7$           & 1  & -2 & -1 & 3  & -1 & 3  & -1 & 2  & 1  & 0 & 0 & 0 & 1 & -1 & 0 & -1 \\
\hline $\Sigma_7'$          & 1  & -2 & -1 & 3  & 1  & -3 & -1 & 2  & 1  & 0 & 0 & 0 & -1 & -1 & 0 & 1 \\
\hline $\Sigma_{14}$        & 14 & 5  & -2 & 2  & -2 & 2  & 2  & 1  & -1 & -1 & -1 & -1 & 0 & 0 & 0 & 1 \\
\hline $\Sigma_{14}'$       & 14 & -4 & 6  & -2 & 0  & 0  & 2  & 0  & 2  & 0 & 0 & -2 & 0 & 0 & 0 & 0 \\
\hline $\Sigma_{14}''$      & 14 & 5  & -2 & 2  & 2  & -2 & 2  & 1  & -1 & 1 & 1 & -1 & 0 & 0 & 0 & -1 \\
\hline $\Sigma_{21}$        & 21 & 3  & 5  & 1  & 3  & -1 & 1  & -1 & 0  & -1 & -1 & 1 & 1 & -1 & 0 & 0 \\
\hline $\Sigma_{21}'$       & 21 & 3  & 5  & 1  & -3 & 1  & 1  & -1 & 0  & 1 & 1 & 1 & -1 & -1 & 0 & 0 \\
\hline $\Sigma_{27}$        & 27 & 0  & 3  & 3  & 3  & 3  & -1 & 0  & 0  & 0 & 0 & 0 & -1 & 1 & -1 & 0 \\
\hline $\Sigma_{27}'$       & 27 & 0  & 3  & 3  & -3 & -3 & -1 & 0  & 0  & 0 & 0 & 0 & 1 & 1 & -1 & 0 \\
\hline $\Sigma_{42}$        & 42 & 6  & 2  & -6 & 0  & 0  & -2 & 2  & 0  & 0 & 0 & 0 & 0 & 0 & 0 & 0 \\
\hline $\Sigma_{56}$        & 56 & 2  & -8 & 0  & 0  & 0  & 0  & -2 & 2  & 0 & 0 & 0 & 0 & 0 & 0 & 0 \\
\hline $\Sigma_{64}$        & 64 & -8 & 0  & 0  & 0  & 0  & 0  & 0  & -2 & 0 & 0 & 0 & 0 & 0 & 1 & 0 \\
\hline
\end{tabular} \\
\caption{Character table for $G_2(2)$.} \label{table:Character_table-IP5}
\end{center}
\end{table}

The group $G_2(2)$ has only two seven-dimensional real representations, both of which are irreducible. Of these two, only $\Sigma_7$ has character values in $[-2,7]$ for all $g \in G_2(2)$, making it the restriction $\varrho_1$ of the seven-dimensional fundamental representation $\rho_1$ of $G_2$ to $G_2(2)$.

\begin{figure}[tb]
\begin{center}
  \includegraphics[width=70mm]{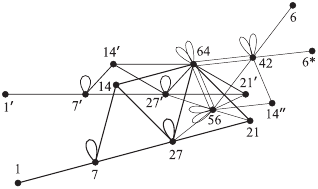}\\
 \caption{The McKay graphs $\mathcal{G}^{\varrho_1}_{G_2(2)}$.} \label{Fig-McKay_Graph-IP5-rho1}
\end{center}
\end{figure}

Thus, up to conjugacy, there is only one embedding of $G_2(2)$ in $G_2$.
The McKay graph for $\varrho_1$ is given in Figure \ref{Fig-McKay_Graph-IP5-rho1}.
This graph is a $\bbZ_2$-orbifold of the McKay graph $\mathcal{G}^{\varrho_1}_{PU(3;3)}$ for $PU(3;3)$.
Note that the McKay graph given in \cite[Figure 2]{he:2003} is not that for $\varrho_1$, but rather for the other irreducible seven-dimensional representation $\Sigma_7'$ of $G_2(2)$.

The decomposition of the Kronecker square of $\varrho_1$ into irreducibles is given by
$$\varrho_1^2 = \mathrm{id} + \varrho_1 + \Sigma_{14} + \Sigma_{27},$$
thus the fourteen-dimensional representation $\varrho_2$ is given by the irreducible representation $\Sigma_{14}$.
We note that the eigenvalues of the representation matrices for $C_4, C_4^{\prime}, C_4^{\prime\prime}$ all square to those for $C_2$, those for $C_6$ square to those for $C_3$, those for $C_6^{\prime}$ square to those for $C_3^{\prime}$, those for $C_8$ square to those for $C_4$, those for $C_8^{\prime}$ square to those for $C_4^{\prime\prime}$, and those for $C_{12}, C_{12}^{\prime}, C_{12}^{\prime\prime}$ all square to those for $C_6$.

Then from (\ref{eqn:moments-subgroupG2}) and Tables \ref{table:Character_table-IP5} and \ref{Table:subgroupsG2-orbits(theta1,theta2)}, we see that
\begin{align*}
\varsigma_{m,n} & = \frac{1}{12096} \Omega^W(0,0) + \frac{56}{12096} \Omega^W(1/3,2/3) + \frac{63+252}{12096} \Omega^W(0,1/2) + \frac{672}{12096} \Omega^W(0,1/3) \\
& \quad + \frac{126+252}{12096} \Omega^W(0,1/4) + \frac{378}{12096} \Omega^W(1/4,1/2) + \frac{504}{12096} \Omega^W(1/6,1/3) \\
& \quad + \frac{1728}{12096} \Omega^W(1/7,3/7) + \frac{1512}{12096} \Omega^W(1/8,1/2) + \frac{1512}{12096} \Omega^W(1/8,3/8) \\
& \quad + \frac{1008+1008+1008}{12096} \Omega^W(1/12,5/12) + \frac{2016}{12096} \Omega^W(1/6,1/2),
\end{align*}
where $\Omega^W(\theta_1,\theta_2)$ is as in Section \ref{sect:II1}.
Thus we obtain:

\begin{Thm}
The joint spectral measure (over $\bbT^2$) for all embeddings of $G_2(2)$ into the fundamental representations of $G_2$ is
\begin{equation}
\begin{split}
\mathrm{d}\varepsilon & = \frac{1}{768\pi^4} J^2 \, \mathrm{d}_8 \times \mathrm{d}_8 + \frac{1}{1344\pi^4} J^2 \, \mathrm{d}_7 \times \mathrm{d}_7 + \frac{1}{1152\pi^4} J^2 \, \mathrm{d}_6 \times \mathrm{d}_6 + \frac{1}{12} \, \mathrm{d}_4 \times \mathrm{d}_4 \\
& \quad + \frac{1}{12} \, \mathrm{d}_3 \times \mathrm{d}_3 - \frac{1}{72} \, \mathrm{d}_2 \times \mathrm{d}_2 - \frac{1}{252} \, \mathrm{d}_1 \times \mathrm{d}_1 + \frac{1}{768\pi^4} J^2 \, \mathrm{d}^{(4)} + \frac{1}{12} \, \mathrm{d}^{(2)} - \frac{1}{24} \, \mathrm{d}^{(1)},
\end{split}
\end{equation}
where $\mathrm{d}_m$ is the uniform measure over $m^{\mathrm{th}}$ roots of unity and $\mathrm{d}^{(k+4)}$ is the uniform measure on the points in $C_k^W$.
\end{Thm}

\subsection{Conclusion and Outlook}

The spectral measures studied here encode the character table of the groups. In the context of quantum subgroups, as discussed in the Introduction and Section \ref{sect:spec_measure-nimrepsG2}, the spectral measure provides the appropriate generalisation of the character table. Whilst the spectral measures for known nimrep graphs associated to the $G_2$ modular invariant partition functions were studied in \cite{evans/pugh:2012i}, the classification of quantum subgroups of $G_2$ is still open, including the classification of all nimrep graphs for $G_2$. And beyond $G_2$, there are many other rank two Lie groups where little is known about their quantum subgroups, including the semisimple, connected, simply-connected, compact group $SU(2) \times SU(2)$, as well as other connected (but not necessarily simply-connected) compact groups, e.g. $U(2)$, $PSU(3)$ or $\bbT \times SO(3)$ \cite{evans/pugh:2012iv}.

The spectral measures for rank two Lie groups and their quantum subgroups are naturally given as measures over $\bbT^2$. The spectral measures presented in this paper for finite subgroups of $G_2$ have been given over $\bbT^2$, however by transporting the measure on a suitable (possibly noncommutative) manifold the picture might simplify. This would be a valuable direction for future research in this area.

\bigskip \bigskip

\begin{footnotesize}
\noindent{\it Acknowledgements.}

The authors thank Robert Griess, Jr. and Nick Gill for their helpful correspondence.
The revision was carried out during EPSRC grant EP/N022432/1.
The second author was supported by the Coleg Cymraeg Cenedlaethol.
\end{footnotesize}

\end{document}